\documentclass{article}

\usepackage{latexsym,calc}
\usepackage[intlimits]{amsmath}
\usepackage{amsthm,amssymb,amscd,wasysym,eufrak}
\usepackage{pst-3dplot}
\usepackage{hyperref}

\renewcommand{\sectionmark}[1]%
  {\markright{\thesection.\ #1}}

\newtheorem{satz}{Theorem}[section]
\newtheorem{lemma}[satz]{Lemma}
\newtheorem{prop}[satz]{Proposition}
\newtheorem{kor}[satz]{Corollary}
\newtheorem{folg}[satz]{Conclusion}

\newcommand{\Satz}[1]{\begin{satz} {\itshape #1}\end{satz}}
\newcommand{\Lemma}[1]{\begin{lemma} {\itshape #1}\end{lemma}}
\newcommand{\Kor}[1]{\begin{kor} {\itshape #1}\end{kor}}

\theoremstyle{definition} 

\newtheorem{bem}[satz]{Note}
\newtheorem{dfn}[satz]{Definition}
\newtheorem{notation}[satz]{Notation}
\newtheorem{example}[satz]{Example}
\newtheorem{examples}[satz]{Examples}

\newcommand{\Bem}[1]{\begin{bem}{#1}\end{bem}}

\newenvironment{beweis}[1][]{\noindent\textit{Proof#1:}%
}

\newcommand{\Beweis}[2][]{\begin{beweis}[#1] {\normalfont #2}{\hfill$\Box$} \end{beweis}}
\newcommand{\Beweisenum}[2][]{\begin{beweis}[#1] {\normalfont #2} \end{beweis}}

\theoremstyle{definition} 

\newtheorem*{convention}{Convention}
\newcommand{\Konv}[1]{\begin{convention}{\normalfont #1}\end{convention}}

\newtheorem*{conclusion}{Conclusion}

\newcommand{\R}{\ensuremath{\mathbb{R}}}

\newcommand{\calC}{\ensuremath{\mathcal{C}}}

\newcommand{\calL}{\ensuremath{\mathcal{L}}}

\newcommand{\oldA}{\ensuremath{\mathfrak{A}}}
\newcommand{\oldB}{\ensuremath{\mathfrak{B}}}
\newcommand{\oldC}{\ensuremath{\mathfrak{C}}}
\newcommand{\olda}{\ensuremath{\mathfrak{a}}}
\newcommand{\oldb}{\ensuremath{\mathfrak{b}}}
\newcommand{\oldc}{\ensuremath{\mathfrak{c}}}

\newcommand{\spann}{{\mbox{span}}}

\newcommand{\sign}{{\mbox{sign}}}
\newcommand{\arcosh}{{\ensuremath{\mathrm{arcosh}}}}

\newcommand{\Hom}{{\ensuremath{\mathrm{Hom}}}}

\newcommand{\dd}{{\makebox[\width]{\white$\mathbb{P}$}\pspicture(0,0)(0,0)\rput[lt]{180}{$\mathbb{P}$}\endpspicture}}

\renewcommand{\phi}{\varphi}
\renewcommand{\epsilon}{\varepsilon}

\newcommand{\llVert}{{\ensuremath{\lvert\!\lvert\!\lvert}}}

\newcommand{\leftllVert}{\ensuremath{\left\lvert\!\left\lvert\!\left\lvert}}

\newcommand{\rrVert}{{\ensuremath{\rvert\!\rvert\!\rvert}}}

\newcommand{\rightrrVert}{\ensuremath{\right\rvert\!\right\rvert\!\right\rvert}}

\newcommand{\llangle}{{\ensuremath{\langle\!\langle}}}

\newcommand{\leftllangle}{\ensuremath{\left\langle\!\left\langle}}

\newcommand{\rrangle}{{\ensuremath{\rangle\!\rangle}}}

\newcommand{\rightrrangle}{\ensuremath{\right\rangle\!\right\rangle}}

\newcommand{\pperp}{{\ensuremath{\mathop{\perp\:\!\!\!\!\!\perp}}}}
\newcommand{\JJ}{{\ensuremath{\mathbb{J}}}}
\renewcommand{\SS}{{\ensuremath{\mathbb{HS}}}}
\newcommand{\ooverline}[1]{{\ensuremath{{\overline{#1}}}}}

\renewcommand{\mid}{\ensuremath{\,|\,}} 

\newcommand{\DEF}[1]{\emph{#1}}

\newcommand{\spacevec}[3]{%
\ensuremath{\begin{pmatrix} #1 \\ #2 \\ #3 \end{pmatrix}}}

\newrgbcolor{grue}{0 .5 .5}
\newcmykcolor{lcyan}{0.3 0 0 0}
\newcmykcolor{purple}{0.2 0.5 0 0.3}
\newrgbcolor{lightred}{1 0.7 0.7}
\newrgbcolor{lviolet}{1 0.5 1}
\newrgbcolor{orange}{1 0.5 0}

\title{{\bf Duality between Hyperbolic and de Sitter Geometry}}
\author{Immanuel Asmus,\\
        University of Potsdam}

\date{\today}

\begin{document}
\sloppy

\maketitle

\begin{abstract}
In this paper we describe trigonometry on the de Sitter surface.
For that a characterization of geodesics is given, leading to various types of triangles.
We define lengths and angles of these.
Then, transferring the concept of polar triangles from spherical geometry into the Minkowski space,
we relate hyperbolic with de Sitter triangles
such that the proof of the hyperbolic law of cosines for angles
becomes much clearer and easier than it is traditionally.
Furthermore, polar triangles turn out to be a powerful tool
for describing de Sitter trigonometry.
\end{abstract}

\tableofcontents

\section*{Notation} \addcontentsline{toc}{section}{Notation}

Throughout this paper we use the following notations, mostly without further explanation in the main body:

\begin{tabular}{p{0.05\textwidth}p{0.85\textwidth}}
 $\llangle.,.\rrangle$ & the Minkowski product on $\R^3$, i.e.\ the bilinear form given by the matrix $\JJ:=(-e_1,e_2,e_3)$, where $\{e_1,e_2,e_3\}$ is the standard basis of $R^3$ \\
 $\llVert.\rrVert$ & the Minkowski (pseudo-)norm: $\llVert x \rrVert=\sqrt{\lvert\llangle x,x\rrangle\rvert}$ \\
 $\calC$ & the light cone, i.e.\ the solution set of $\llVert x\rrVert=0$ \\
 $\SS^2$ & the solution set of $\llVert x \rrVert=1$, consisting of the de Sitter surface $S^{1,1}$ and two copies of the hyperbolic surface, which we denote by $H^2$ (the part containing $e_1$) and $(-H^2)$, respectively \\
 $\hat{x}$ & $x$ divided by its Minkowski norm. Obviously, $x$ must not be in $\calC$. \\
 $\pperp$ & Minkowski or Lorentz orthogonal: $x\pperp y:\Leftrightarrow\llangle x,y\rrangle=0$ \\
 $\calL(3)$ & the Lorentz group. Its elements are the Lorentz transformations. \\
 $d_H$ & {\raggedright the hyperbolic distance, a metric on $H^2$:\\$d_H(x,y):=\arcosh(-\llangle x,y \rrangle)$} \\
 $d_H'$ & {\raggedright the antipodal hyperbolic distance, a metric on $(-H^2)$:\\$d_H'(x,y):=d_H(-x,-y)$} \\
 $\dd$ & the proper de Sitter distance, a pseudometric on $S^{1,1}$ (will be defined later) \\
 $d_\SS$ & the generalized de Sitter distance, which equals $d_H$, $d_H'$ or $\dd$ depending on the surface that the points are located on \\
 $\ooverline{ab}$ & the generalized de Sitter segment (defined later)
\end{tabular}

\section{Introduction}

There exists a vast variety of books dealing with hyperbolic trigonometry.
They give one or two laws of cosines and the law of sines,
and almost each of these books provides another proof for those theorems.
The literature given in the references section shows some typical ways to prove the hyperbolic trigonometric rules.

The most elementary proofs may be found in Wilson \cite{Wil} and Anderson \cite{And}.
Wilson constructs a triangle with given side lengths in the hyperboloid model,
where the first point equals $e_1$, the second one is located in the $e_1$-$e_2$-plane,
and the third one has a positive third coordinate.
A simple computation yields the law of cosines for sides and the law of sines.
By congruence, these results can be generalized to arbitrary triangles.
The law of cosines for angles is not mentioned.
Anderson goes mainly the same way, except that he uses the Poincar\'e disc model,
where the constructed triangle has one point in the origin,
the second one on the positive real axis, and the third one has positive imaginary part.
He directly derives the law of cosines for sides.
For both of the other trigonometric laws he needs a purely algebraic but not obvious computation.

Iversen \cite{Ive} and Ungar \cite{Ung} prove the trigonometric rules in a more direct and simple way,
but use rather abstract models.
Iversen describes hyperbolic geometry in the $sl_2$ model,
where each point of $H^2$ is given by a $2\times2$ matrix with vanishing trace, determinant equal to 1
and positive lower left component.
Ungar uses the gyrovector space, where even the trigonometric rules themselves take an almost unrecognizable form.

Finally, Thurston \cite{Thu} works in the hyperboloid model.
He uses the fact that a hyperbolic triangle $\{v_1,v_2,v_3\}$ forms a basis of $\R^3$,
computes the dual basis $\{w_1,w_2,w_3\}$
and then shows that the matrix $\bigl(\llangle w_i,w_j\rrangle\bigr)_{i,j}$ is the inverse of
$\bigl(\llangle v_i,v_j\rrangle\bigr)_{i,j}$.
From this he obtains the law of cosines for sides.
He repeats the same with the $v_i$ lying on the de Sitter surface.
However, Thurston does not use these vectors to describe de Sitter geometry, but to describe hyperbolic geodesics
(or hyperplanes, in higher dimensions).
To then obtain the law of cosines for angles, he needs to take enhanced care of the occuring signs.
The law of sines is derived by algebraic transformations from the law of cosines for sides,
firstly only for right triangles.
Dividing an arbitrary triangle into two right triangles by one of its altitudes
and applying the known relation to these right triangles yields the general law of sines.

\vspace{.5cm}
In this work we only use the hyperboloid model, in which we are able to
provide simple and illustrative proofs for all of the trigonometric laws for hyperbolic geometry.
Our proof for the hyperbolic law of cosines for sides is similar to the one given by Iversen.
From the paremetrization of hyperbolic geodesics,
we can compute unit tangent vectors at one vertex of a triangle, pointing in the direction of another vertex:
$$ X_{AB}=\frac{B+\llangle A,B\rrangle A}{\llVert A\times B\rrVert}. $$
We find the size of an angle at a vertex, according to its definition,
by applying the Minkowski product to the respective tangent vectors:
$$ \cos(\alpha)=\llVert X_{AB},X_{AC}\rrVert. $$
This yields directly the hyperbolic law of cosines for sides
$$ \cosh(a)=\cosh(b)\cosh(c)-\cos(\alpha)\sinh(b)\sinh(c). $$

The same can be done for a certain kind of triangles on the de Sitter surface,
which we call non-contractible spatiolateral triangles.
We obtain for these triangles the following law of cosines for sides:
$$ \cos(a)=\cos(b)\cos(c)-\cosh(\alpha)\sin(b)\sin(c). $$

By the duality between the hyperbolic plane and the de Sitter surface
we get the relation $\alpha=\pi-a'$ between the angles of a hyperbolic triangle $\{A,B,C\}$
and the sides of the associated polar triangle $\{A',B',C'\}$,
which corresponds to the dual basis used in Thurston's work, save that they are normalized:
$$ A':=\det(A,B,C)\frac{B\times C}{\llVert B\times C\rrVert}. $$
Conversely, we relate the sides of the hyperbolic triangle with the angles of the polar triangle and get $a=\alpha'$.
Plugging in these relations into the law of cosines for sides in non-contractible spatiolateral triangles,
we get the hyperbolic law of cosines for angles,
$$ \cos(\alpha)=-\cos(\beta)\cos(\gamma)+\cosh(a)\sin(\beta)\sin(\gamma). $$

To obtain the hyperbolic law of sines, we again use the duality to see that
$$ \sin(\alpha)=\sin(a')=\frac{\lvert\det(A,B,C)\rvert}{\llVert A\times B\rrVert\cdot\llVert A\times C\rrVert}. $$
Dividing by $\sinh(a)=\llVert B\times C\rrVert$ gives a term that is symmetric in $A$, $B$, and $C$,
which proves the theorem.

On the way to these results, we get all these trigonometric rules for three types of de Sitter triangles,
called contractible spatiolateral, non-contractible spatiolateral, and tempolateral triangles.
Furthermore, we are able to make statements about the sum of the side lengths in spatiolateral triangles,
and we show that in every contractible spatiolateral triangle there exists one side
longer than the sum of the other sides.

\section*{Acknowledgement}

First of all I would like to thank my supervisor Christian B\"ar.
He has taught me the fundamentals of differential geometry and encouraged me to write this work.
Thanks are also due to my room mates Florian Hanisch and Christian Becker
who assisted me with the translation.

\section{Fundamentals of Minkowski Geometry}

When we talk about Minkowski space, we think of $\R^n$ with a non-degenerate, symmetric bilinear form $\llangle.,.\rrangle$, where the matrix representing $\llangle.,.\rrangle$ has the eigenvalues 1 with multiplicity $(n-1)$, and $-1$ with multiplicity 1.
In this paper we deal with three-dimensional Minkowski space exclusively,
because when observing triangles on a higher than two-dimensional hyperbolic or de Sitter surface,
each such triangle lies completely in a two-dimensional submanifold
that is either hyperbolic, de Sitter, or spherical.

In the two-dimensional case, thus, the bilinear form is represented by $(-e_1,e_2,e_3)$, as given in the ``notation'' section.
Because of its indefiniteness, the Minkowski product gives a partition of $\R^3$, according to the sign of the associated quadratic form.
A vector $x\in\R^3$ is called
\begin{itemize}
 \item \emph{timelike,} if $\llangle x,x\rrangle<0$;
 \item \emph{lightlike,} if $\llangle x,x\rrangle=0$ and $x\neq0$; and
 \item \emph{spacelike,} if $\llangle x,x\rrangle>0$ or $x=0$, respectively.
\end{itemize}
Timelike and lightlike vectors are embraced by the term \emph{``causal''} vectors.
The lightlike vectors together with 0 form the \emph{light cone} $\calC$.

Obviously, lightlike vectors cannot be normalized.
The set of normalized vectors, denoted by $\SS^2$, consists of three components, as figure~\ref{C,SS2} shows.

\begin{figure}
\begin{center}
{\psset{unit=0.65cm}
   \pspicture(-3,-3.5)(3,3.7) 
    \psset{Beta=15}
    \parametricplotThreeD[plotstyle=curve,hiddenLine=true,linecolor=cyan,linewidth=0.2pt](135,315)(-2,2){%
     t cos u dup mul 1 add sqrt mul t sin u dup mul 1 add sqrt mul u}
    \parametricplotThreeD[plotstyle=curve,hiddenLine=true,linecolor=cyan,linewidth=0.2pt](-2,2)(135,315){%
     u cos t dup mul 1 add sqrt mul u sin t dup mul 1 add sqrt mul t}

    \parametricplotThreeD[plotstyle=curve,hiddenLine=true,linecolor=yellow,linewidth=0.2pt](135,315)(0,2.25){%
     t cos u mul t sin u mul u}
    \parametricplotThreeD[plotstyle=curve,hiddenLine=true,linecolor=yellow,linewidth=0.2pt](0,2.25)(135,315){%
     u cos t mul u sin t mul t}
    \parametricplotThreeD[plotstyle=curve,hiddenLine=true,linecolor=yellow,linewidth=0.2pt](-45,135)(-2.25,0){%
     t cos u mul t sin u mul u}
    \parametricplotThreeD[plotstyle=curve,hiddenLine=true,linecolor=yellow,linewidth=0.2pt](-2.25,0)(-45,135){%
     u cos t mul u sin t mul t}

    \parametricplotThreeD[plotstyle=curve,hiddenLine=true,linecolor=lightred](0,360)(1,2.5){%
     t cos u dup mul 1 neg add sqrt mul t sin u dup mul 1 neg add sqrt mul u}
    \parametricplotThreeD[plotstyle=curve,hiddenLine=true,linecolor=lightred](0,360)(-1,-2.5){%
     t cos u dup mul 1 neg add sqrt mul t sin u dup mul 1 neg add sqrt mul u}
    \parametricplotThreeD[plotstyle=curve,hiddenLine=true,linecolor=lightred](1,2.5)(0,360){%
     u cos t dup mul 1 neg add sqrt mul u sin t dup mul 1 neg add sqrt mul t}
    \parametricplotThreeD[plotstyle=curve,hiddenLine=true,linecolor=lightred](-1,-2.5)(0,360){%
     u cos t dup mul 1 neg add sqrt mul u sin t dup mul 1 neg add sqrt mul t}

    \parametricplotThreeD[plotstyle=curve,hiddenLine=true,linecolor=yellow,linewidth=0.2pt](-45,135)(-2.25,2.25){%
     t cos u mul t sin u mul u}
    \parametricplotThreeD[plotstyle=curve,hiddenLine=true,linecolor=yellow,linewidth=0.2pt](-2.25,2.25)(-45,135){%
     u cos t mul u sin t mul t}
    \parametricplotThreeD[plotstyle=curve,hiddenLine=true,linecolor=yellow,linewidth=0.2pt](135,315)(-2.25,0){%
     t cos u mul t sin u mul u}
    \parametricplotThreeD[plotstyle=curve,hiddenLine=true,linecolor=yellow,linewidth=0.2pt](-2.25,0)(135,315){%
     u cos t mul u sin t mul t}

    \parametricplotThreeD[plotstyle=curve,hiddenLine=true,linecolor=cyan,linewidth=0.2pt](-2,2)(-45,135){%
     u cos t dup mul 1 add sqrt mul u sin t dup mul 1 add sqrt mul t}
    \parametricplotThreeD[plotstyle=curve,hiddenLine=true,linecolor=cyan,linewidth=0.2pt](-45,135)(-2,2){%
     t cos u dup mul 1 add sqrt mul t sin u dup mul 1 add sqrt mul u}

    \pstThreeDCoor[nameX=$x_2$,nameY=$x_3$,nameZ=$x_1$,zMin=-3,zMax=3,xMin=-4,yMin=-4,linewidth=0.2pt,linecolor=black]
    \pstThreeDPut(-2,2,2.7){\red $H^2$} \pstThreeDPut(-2,2,2.2){\orange $\calC$~~} \pstThreeDPut(-2,2,1.7){\blue $S^{1,1}$}
   \endpspicture
}
\caption{The light cone $\calC$ and the surfaces $\SS^2$ of normalized vectors} \label{C,SS2}
\end{center}
\end{figure}

The outer surface is the de Sitter surface $S^{1,1}$, which is a Lorentz manifold,
while the inner ones are two copies of the hyperbolic surface, which is a Riemannian manifold.
According to the sign of the first component, the copies are denoted by $H^2$ and $(-H^2)$, respectively.

We call two vectors $x$ and $y$ \emph{Lorentz} or \emph{Minkowski orthogonal,} if $\llangle x,y \rrangle=0$.
Euclidean orthogonality is neither a necessary nor a sufficient condition for Lorentz orthogonality, but they may go together.
For example, $(1,1,0)^t$ is Minkowski orthogonal to itself for it is lightlike,
whereas $(1,1,0)^t$ and $(-1,1,0)^t$ are Euclidean, but not Lorentz orthogonal.
Vectors $e_1$ and $e_2$ are orthogonal on both counts.

\Bem{\label{bem:lichtlicht} Two lightlike vectors are Minkowski orthogonal if and only if they are linearly dependent.
To obtain this, we write
 $$x=\spacevec{x_1}{x_1\cos\alpha}{x_1\sin\alpha} \quad \mbox{and} \quad
 y=\spacevec{y_1}{y_1\cos\beta}{y_1\sin\beta}.$$
Computing
 $$0=\llangle{}x,y\rrangle=x_1y_1(\cos\alpha\cos\beta+\sin\alpha\sin\beta-1)$$
leads to
 $$\alpha=\beta+2k\pi,$$
and thus $\sin\alpha=\sin\beta$ and $\cos\alpha=\cos\beta$ hold, showing linear dependence.}

\Lemma{\label{satz:exist_MOB} Let $U\subset\R^3$ be a vector subspace.
 Then, there exists a Lorentz orthogonal basis of $U$, i.e.~a basis $\{b_i\}$ which satisfies $b_i\pperp b_j$ whenever $i\neq j$.}

\Beweis{Since $\llangle.,.\rrangle|_U$ is symmetric, the matrix representing this bilinear form can be diagonalized.
This is equivalent to the statement of the lemma.}

\vspace{.5cm}
We call a basis $\{b_i\}$ \emph{Lorentz} or \emph{Minkowski orthonormal,} if it is Lorentz orthogonal and
 $$ \llangle{}b_i,b_i\rrangle=\begin{cases}-1,&i=1\\1,&otherwise\end{cases}$$
holds.
The Lorentz group
 \[\calL(3):=\left\{\Phi\in\Hom(\R^3,\R^3)\mid\llangle\Phi(x),\Phi(y)\rrangle=\llangle x,y \rrangle\;\forall x,y\in\R^3\right\}\]
is therefore equivalently characterized by consisting of matrices $(b_1,b_2,b_3)$ with $\{b_1,b_2,b_3\}$ being a Minkowski orthonormal basis of $\R^3$.

\Lemma{\label{folg:Basis_zeit_2raum} Each Lorentz orthogonal basis of $\R^3$ consists of one timelike vector and two spacelike vectors.}

\Beweis{This again is pure linear algebra, namely Sylvester's law of inertia applied to $\llangle.,.\rrangle$.}

\vspace{.5cm}
 It can be easily proven that any nonzero vector Lorentz orthogonal to a timelike one must itself be spacelike (cf.~Naber \cite{Nab}),
 so together with note~\ref{bem:lichtlicht} one can conclude that every two-dimensional subspace of $\R^3$
 contains a nonzero spacelike vector.
 This finding allows us to classify those subspaces by the type of the second vector in a basis.

\dfn{Let $U$ be a two-dimensional vector subspace of $R^3$.
 Let further $\{b_1,b_2\}$ be a Minkowski orthogonal basis of $U$, with $b_1$ being spacelike.
 We call $U$
 \begin{enumerate}
  \item \DEF{spacelike,} if $b_2$ is spacelike as well;
  \item \DEF{lightlike,} if $b_2$ is lightlike; and
  \item \DEF{timelike,} if $b_2$ is timelike.
 \end{enumerate}}

This classification is independent of the choice of the basis,
and the property of a plane to be spacelike, timelike, or lightlike, respectively, does not change under Lorentz transformations.

\begin{figure}
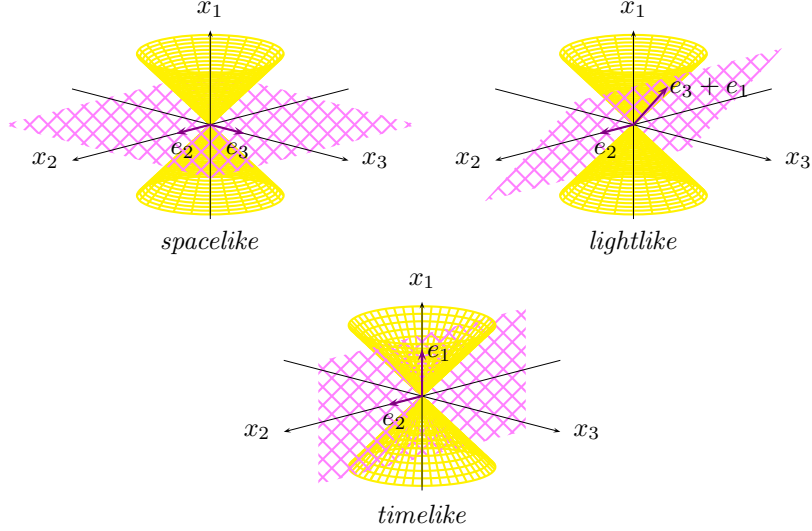

 \begin{center}
  {\psset{unit=0.65cm}
    \pspicture(-5,-3.5)(3,3.25) 
     \psset{Beta=15}
     \parametricplotThreeD[plotstyle=curve,hiddenLine=true,linecolor=yellow](0,360)(-1.5,0){%
      t cos u mul t sin u mul u}
     \parametricplotThreeD[plotstyle=curve,hiddenLine=true,linecolor=yellow](-1.5,0)(0,360){%
      u cos t mul u sin t mul t}
     \pstThreeDSquare[linestyle=none,fillstyle=crosshatch,hatchcolor=lviolet](-3,-3,0)(6,0,0)(0,6,0)
     \parametricplotThreeD[plotstyle=curve,hiddenLine=true,linecolor=yellow](0,360)(0,1.5){%
      t cos u mul t sin u mul u}
     \parametricplotThreeD[plotstyle=curve,hiddenLine=true,linecolor=yellow](0,1.5)(0,360){%
      u cos t mul u sin t mul t}
     \pstThreeDCoor[nameX=$x_2$,nameY=$x_3$,nameZ=$x_1$,zMin=-2,zMax=2,xMin=-4,yMin=-4,linewidth=0.2pt,linecolor=black]
     \pstThreeDLine[linecolor=violet,arrows=<->](1,0,0)(0,0,0)(0,1,0)
     \pstThreeDPut(1.3,0.5,-0.2){$e_2$} \pstThreeDPut(0.5,1.3,-0.2){$e_3$}
     \pstThreeDPut(0,0,-2.5){\emph{spacelike}}
   \endpspicture\hspace{\stretch{1}}
   \pspicture(-5,-3.5)(5,3.25) 
    \psset{Beta=15}
    \parametricplotThreeD[plotstyle=curve,hiddenLine=true,linecolor=yellow](0,360)(0,1.5){%
     t cos u mul t sin u mul u}
    \parametricplotThreeD[plotstyle=curve,hiddenLine=true,linecolor=yellow](0,1.5)(0,360){%
     u cos t mul u sin t mul t}
    \pstThreeDSquare[linestyle=none,fillstyle=crosshatch,hatchcolor=lviolet](-3,-1.3,-1.3)(0,2.6,2.6)(6,0,0)
    \parametricplotThreeD[plotstyle=curve,hiddenLine=true,linecolor=yellow](0,360)(-1.5,0){%
     t cos u mul t sin u mul u}
    \parametricplotThreeD[plotstyle=curve,hiddenLine=true,linecolor=yellow](-1.5,0)(0,360){%
     u cos t mul u sin t mul t}
    \pstThreeDCoor[nameX=$x_2$,nameY=$x_3$,nameZ=$x_1$,zMin=-2,zMax=2,xMin=-4,yMin=-4,linewidth=0.2pt,linecolor=black]
    \pstThreeDLine[linecolor=violet,arrows=<->](1,0,0)(0,0,0)(0,1,1)
    \pstThreeDPut(1.3,0.5,-0.2){$e_2$} \pstThreeDPut(0,2.2,1.2){$e_3+e_1$}
    \pstThreeDPut(0,0,-2.5){\emph{lightlike}}
   \endpspicture

   \pspicture(-3,-2)(3,2) 
    \psset{Beta=15}
    \parametricplotThreeD[plotstyle=curve,hiddenLine=true,linecolor=yellow](180,360)(-1.5,1.5){%
     t cos u abs mul t sin u abs mul u}
    \parametricplotThreeD[plotstyle=curve,hiddenLine=true,linecolor=yellow](-1.5,1.5)(180,360){%
     u cos t abs mul u sin t abs mul t}
    \pstThreeDSquare[linestyle=none,fillstyle=crosshatch,hatchcolor=lviolet](-3,0,-1.3)(0,0,2.6)(6,0,0)
    \parametricplotThreeD[plotstyle=curve,hiddenLine=true,linecolor=yellow](0,180)(-1.5,1.5){%
     t cos u abs mul t sin u abs mul u}
    \parametricplotThreeD[plotstyle=curve,hiddenLine=true,linecolor=yellow](-1.5,1.5)(0,180){%
     u cos t abs mul u sin t abs mul t}
    \pstThreeDCoor[nameX=$x_2$,nameY=$x_3$,nameZ=$x_1$,zMin=-2,zMax=2,xMin=-4,yMin=-4,linewidth=0.2pt,linecolor=black]
    \pstThreeDLine[linecolor=violet,arrows=<->](1,0,0)(0,0,0)(0,0,1)
    \pstThreeDPut(1.3,0.5,-0.2){$e_2$} \pstThreeDPut(0,0.5,1){$e_1$}
    \pstThreeDPut(0,0,-2.5){\emph{timelike}}
   \endpspicture
  }
 \end{center}
 \caption{Examples of spacelike, lightlike, and timelike planes}
\end{figure}

\vspace{.5cm}
Geodesics on $\SS^2$ (so-called \emph{generalized de Sitter geodesics}) can be obtained from the intersection of
a two-dimensional vector subspace of $\R^3$ with $\SS^2$:
The intersection of a spacelike or lightlike plane with $S^{1,1}$ will be a great ellipse or a pair of parallel lines, respectively.
Such planes do not intersect with $(\pm H^2)$.
Timelike planes, however, intersect with each component of $\SS^2$, giving two great hyperbolas (i.e.\ four great hyperbola branches).

Geodesics obtained in this way are the usual geodesics according to the Riemannian (for $H^2$)
or semi-Riemannian metric (for $S^{1,1}$).
Geodesics on $H^2$ we call \emph{hyperbolic,} whereas geodesics on $(-H^2)$ will be referred to as being \emph{antipodal hyperbolic.}
Geodesics on $S^{1,1}$ are called \emph{(proper) de Sitter geodesics,}
but we should keep in mind, that they can be of three different kinds: either ellipses, straight lines, or hyperbola branches.

\begin{figure}
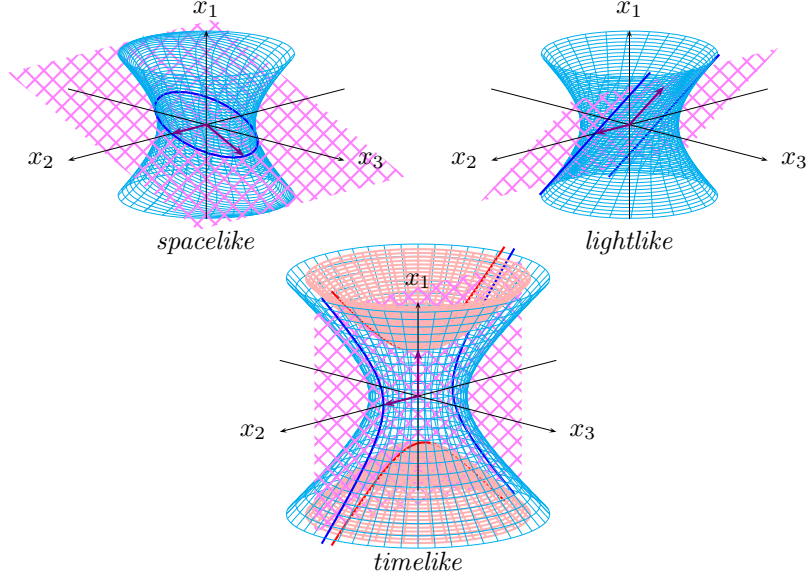

 \begin{center}
  {\psset{unit=0.65cm}
   \pspicture(-5,-3.5)(3,3.25) 
    \psset{Beta=15}
    \parametricplotThreeD[plotstyle=curve,hiddenLine=true,linecolor=cyan,linewidth=0.2pt](0,360)(-1.5,0){%
     t cos u dup mul 1 add sqrt mul
     t sin u dup mul 1 add sqrt mul u 1.5 add u 1.5 neg add mul dup mul 50.625 div 1 add mul u u 1.5 add u 1.5 neg add mul dup mul 50.625 div 1 add dup mul 1 neg add sqrt mul neg add
     t sin u dup mul 1 add sqrt mul neg u 1.5 add u 1.5 neg add mul dup mul 50.625 div 1 add dup mul 1 neg add sqrt mul u u 1.5 add u 1.5 neg add mul dup mul 50.625 div 1 add mul add}
    \parametricplotThreeD[plotstyle=curve,hiddenLine=true,linecolor=cyan,linewidth=0.2pt](-1.5,0)(0,360){%
     u cos t dup mul 1 add sqrt mul
     u sin t dup mul 1 add sqrt mul t 1.5 add t 1.5 neg add mul dup mul 50.625 div 1 add mul t t 1.5 add t 1.5 neg add mul dup mul 50.625 div 1 add dup mul 1 neg add sqrt mul neg add
     u sin t dup mul 1 add sqrt mul neg t 1.5 add t 1.5 neg add mul dup mul 50.625 div 1 add dup mul 1 neg add sqrt mul t t 1.5 add t 1.5 neg add mul dup mul 50.625 div 1 add mul add}
    \pstThreeDSquare[linestyle=none,fillstyle=crosshatch,hatchcolor=lviolet](-3,-2.75,1.1456)(6,0,0)(0,5.5,-2.2913)
    \parametricplotThreeD[plotstyle=curve,hiddenLine=true,linecolor=cyan,linewidth=0.2pt](0,360)(0,1.5){%
     t cos u dup mul 1 add sqrt mul
     t sin u dup mul 1 add sqrt mul u 1.5 add u 1.5 neg add mul dup mul 50.625 div 1 add mul u u 1.5 add u 1.5 neg add mul dup mul 50.625 div 1 add dup mul 1 neg add sqrt mul neg add
     t sin u dup mul 1 add sqrt mul neg u 1.5 add u 1.5 neg add mul dup mul 50.625 div 1 add dup mul 1 neg add sqrt mul u u 1.5 add u 1.5 neg add mul dup mul 50.625 div 1 add mul add}
    \pstThreeDEllipse[linecolor=blue](0,0,0)(1,0,0)(0,1.1,-0.4583)
    \parametricplotThreeD[plotstyle=curve,hiddenLine=true,linecolor=cyan,linewidth=0.2pt](0,1.5)(0,360){%
     u cos t dup mul 1 add sqrt mul
     u sin t dup mul 1 add sqrt mul t 1.5 add t 1.5 neg add mul dup mul 50.625 div 1 add mul t t 1.5 add t 1.5 neg add mul dup mul 50.625 div 1 add dup mul 1 neg add sqrt mul neg add
     u sin t dup mul 1 add sqrt mul neg t 1.5 add t 1.5 neg add mul dup mul 50.625 div 1 add dup mul 1 neg add sqrt mul t t 1.5 add t 1.5 neg add mul dup mul 50.625 div 1 add mul add}
    \pstThreeDCoor[nameX=$x_2$,nameY=$x_3$,nameZ=$x_1$,zMin=-2,zMax=2,xMin=-4,yMin=-4,linewidth=0.2pt,linecolor=black]
    \pstThreeDLine[linecolor=violet,arrows=<->](1,0,0)(0,0,0)(0,1.1,-0.4583)
    \pstThreeDPut(0,0,-2.5){\emph{spacelike}}
   \endpspicture\hspace{\stretch{1}}
   \pspicture(-5,-3.5)(5,3.25) 
    \psset{Beta=15}
    \parametricplotThreeD[plotstyle=curve,hiddenLine=true,linecolor=cyan,linewidth=0.2pt](90,270)(0,1.5){%
     t sin 1 u dup mul 1 add sqrt div arcsin 90 neg add mul t add sin u dup mul 1 add sqrt mul t sin 1 u dup mul 1 add sqrt div arcsin 90 neg add mul t add cos u dup mul 1 add sqrt mul u}
    \parametricplotThreeD[plotstyle=curve,hiddenLine=true,linecolor=cyan,linewidth=0.2pt](0,1.5)(90,270){%
     u sin 1 t dup mul 1 add sqrt div arcsin 90 neg add mul u add sin t dup mul 1 add sqrt mul u sin 1 t dup mul 1 add sqrt div arcsin 90 neg add mul u add cos t dup mul 1 add sqrt mul t}
    \parametricplotThreeD[plotstyle=curve,hiddenLine=true,linecolor=cyan,linewidth=0.2pt](-90,90)(-1.5,0){%
     t sin 1 u dup mul 1 add sqrt div arcsin 90 neg add mul t add sin u dup mul 1 add sqrt mul t sin 1 u dup mul 1 add sqrt div arcsin 90 neg add mul t add cos u dup mul 1 add sqrt mul neg u}
    \parametricplotThreeD[plotstyle=curve,hiddenLine=true,linecolor=cyan,linewidth=0.2pt](-1.5,0)(-90,90){%
     u sin 1 t dup mul 1 add sqrt div arcsin 90 neg add mul u add sin t dup mul 1 add sqrt mul u sin 1 t dup mul 1 add sqrt div arcsin 90 neg add mul u add cos t dup mul 1 add sqrt mul neg t}
    \pstThreeDSquare[linestyle=none,fillstyle=crosshatch,hatchcolor=lviolet](-3,-1.3,-1.3)(0,2.6,2.6)(6,0,0)
    \parametricplotThreeD[plotstyle=curve,hiddenLine=true,linecolor=cyan,linewidth=0.2pt](0,1.5)(-90,90){%
     u sin 1 t dup mul 1 add sqrt div arcsin 90 neg add mul u add sin t dup mul 1 add sqrt mul u sin 1 t dup mul 1 add sqrt div arcsin 90 neg add mul u add cos t dup mul 1 add sqrt mul t}
    \parametricplotThreeD[plotstyle=curve,hiddenLine=true,linecolor=cyan,linewidth=0.2pt](-1.5,0)(90,270){%
     u sin 1 t dup mul 1 add sqrt div arcsin 90 neg add mul u add sin t dup mul 1 add sqrt mul u sin 1 t dup mul 1 add sqrt div arcsin 90 neg add mul u add cos t dup mul 1 add sqrt mul neg t}
    \pstThreeDLine[linecolor=blue](-1,-1.6,-1.6)(-1,1.6,1.6)
    \parametricplotThreeD[plotstyle=curve,hiddenLine=true,linecolor=cyan,linewidth=0.2pt](-90,90)(0,1.5){%
     t sin 1 u dup mul 1 add sqrt div arcsin 90 neg add mul t add sin u dup mul 1 add sqrt mul t sin 1 u dup mul 1 add sqrt div arcsin 90 neg add mul t add cos u dup mul 1 add sqrt mul u}
    \parametricplotThreeD[plotstyle=curve,hiddenLine=true,linecolor=cyan,linewidth=0.2pt](90,270)(-1.5,0){%
     t sin 1 u dup mul 1 add sqrt div arcsin 90 neg add mul t add sin u dup mul 1 add sqrt mul t sin 1 u dup mul 1 add sqrt div arcsin 90 neg add mul t add cos u dup mul 1 add sqrt mul neg u}
    \pstThreeDLine[linecolor=blue](1,-1.6,-1.6)(1,1.6,1.6)
    \pstThreeDCoor[nameX=$x_2$,nameY=$x_3$,nameZ=$x_1$,zMin=-2,zMax=2,xMin=-4,yMin=-4,linewidth=0.2pt,linecolor=black]
    \pstThreeDLine[linecolor=violet,arrows=<->](1,0,0)(0,0,0)(0,1,1)
    \pstThreeDPut(0,0,-2.5){\emph{lightlike}}
   \endpspicture

   \pspicture(-3,-3)(3,2) 
    \psset{Beta=15}
    \parametricplotThreeD[plotstyle=curve,hiddenLine=true,linecolor=cyan,linewidth=0.2pt](180,315)(-2.5,2.5){%
     t cos u dup mul 1 add sqrt mul t sin u dup mul 1 add sqrt mul u}
    \parametricplotThreeD[plotstyle=curve,hiddenLine=true,linecolor=cyan,linewidth=0.2pt,yPlotpoints=19](-2.5,2.5)(180,315){%
     u cos t dup mul 1 add sqrt mul u sin t dup mul 1 add sqrt mul t}
    \parametricplotThreeD[plotstyle=curve,hiddenLine=true,linecolor=lightred](180,360)(1,2.5){%
     t cos u dup mul 1 neg add sqrt mul t sin u dup mul 1 neg add sqrt mul u}
    \parametricplotThreeD[plotstyle=curve,hiddenLine=true,linecolor=lightred](180,360)(-1,-2.5){%
     t cos u dup mul 1 neg add sqrt mul t sin u dup mul 1 neg add sqrt mul u}
    \parametricplotThreeD[plotstyle=curve,hiddenLine=true,linecolor=lightred](1,2.5)(180,360){%
     u cos t dup mul 1 neg add sqrt mul u sin t dup mul 1 neg add sqrt mul t}
    \parametricplotThreeD[plotstyle=curve,hiddenLine=true,linecolor=lightred](-1,-2.5)(180,360){%
     u cos t dup mul 1 neg add sqrt mul u sin t dup mul 1 neg add sqrt mul t}
    \parametricplotThreeD[plotstyle=curve,hiddenLine=true,linecolor=cyan,linewidth=0.2pt](315,360)(-2.5,2.5){%
     t cos u dup mul 1 add sqrt mul t sin u dup mul 1 add sqrt mul u}
    \parametricplotThreeD[plotstyle=curve,hiddenLine=true,linecolor=cyan,linewidth=0.2pt,yPlotpoints=7](-2.5,2.5)(315,360){%
     u cos t dup mul 1 add sqrt mul u sin t dup mul 1 add sqrt mul t}

    \pstThreeDSquare[linestyle=none,fillstyle=crosshatch,hatchcolor=lviolet](-3,0,-2.3)(0,0,4.6)(6,0,0)

    \parametricplotThreeD[plotstyle=curve,hiddenLine=true,linecolor=cyan,linewidth=0.2pt,yPlotpoints=7](-2.5,2.5)(135,180){%
     u cos t dup mul 1 add sqrt mul u sin t dup mul 1 add sqrt mul t}
    \parametricplotThreeD[plotstyle=curve,hiddenLine=true,linecolor=cyan,linewidth=0.2pt](135,180)(-2.5,2.5){%
     t cos u dup mul 1 add sqrt mul t sin u dup mul 1 add sqrt mul u}
    \parametricplotThreeD[plotstyle=curve,hiddenLine=true,linecolor=lightred](1,2.5)(0,180){%
     u cos t dup mul 1 neg add sqrt mul u sin t dup mul 1 neg add sqrt mul t}
    \parametricplotThreeD[plotstyle=curve,hiddenLine=true,linecolor=lightred](-1,-2.5)(0,180){%
     u cos t dup mul 1 neg add sqrt mul u sin t dup mul 1 neg add sqrt mul t}
    \parametricplotThreeD[plotstyle=curve,linecolor=red](-2.5,2.5){t 0 t dup mul 1 add sqrt}
    \parametricplotThreeD[plotstyle=curve,linecolor=red](-2.5,2.5){t 0 t dup mul 1 add sqrt neg}
    \parametricplotThreeD[plotstyle=curve,linecolor=blue](-2.6,2.6){t dup mul 1 add sqrt neg 0 t}
    \parametricplotThreeD[plotstyle=curve,hiddenLine=true,linecolor=lightred](0,180)(1,2.5){%
     t cos u dup mul 1 neg add sqrt mul t sin u dup mul 1 neg add sqrt mul u}
    \parametricplotThreeD[plotstyle=curve,hiddenLine=true,linecolor=lightred](0,180)(-1,-2.5){%
     t cos u dup mul 1 neg add sqrt mul t sin u dup mul 1 neg add sqrt mul u}
    \parametricplotThreeD[plotstyle=curve,hiddenLine=true,linecolor=cyan,linewidth=0.2pt,yPlotpoints=19](-2.5,2.5)(0,135){%
     u cos t dup mul 1 add sqrt mul u sin t dup mul 1 add sqrt mul t}
    \parametricplotThreeD[plotstyle=curve,linecolor=blue](-2.6,2.6){t dup mul 1 add sqrt 0 t}
    \parametricplotThreeD[plotstyle=curve,hiddenLine=true,linecolor=cyan,linewidth=0.2pt](0,135)(-2.5,2.5){%
     t cos u dup mul 1 add sqrt mul t sin u dup mul 1 add sqrt mul u}
    \pstThreeDCoor[nameX=$x_2$,nameY=$x_3$,nameZ=$x_1$,zMin=-2,zMax=2,xMin=-4,yMin=-4,linewidth=0.2pt,linecolor=black]
    \pstThreeDLine[linecolor=violet,arrows=<->](1,0,0)(0,0,0)(0,0,1)
    \pstThreeDPut(0,0,-3.5){\emph{timelike}}
   \endpspicture
  }
 \end{center}
 \caption{Intersections of $\SS^2$ with a spacelike, lightlike, and timelike plane, leading to different types of geodesics}
\end{figure}

\vspace{.5cm}
For a metric on $H^2$, we choose the standard Riemannian metric
\[d_H(x,y):=\arcosh\bigl(-\llangle x,y\rrangle\bigr) \quad \forall x,y\in H^2,\]
which is called \emph{hyperbolic distance.}
The same is done the most simple way for $(-H^2)$, leading to the \emph{antipodal hyperbolic distance}
\[d_H'(x,y):=d_H(-x,-y) \quad \forall x,y\in (-H^2).\]
One can easily see, that both distance functions are well-defined and are indeed metrics.

Things are getting more complicated when the de Sitter surface is concerned.
On the hyperbolic plane, every two points are located on some great hyperbola branch,
thus making it sensible to define their distance by using hyperbolic functions.
However, on de Sitter surface, two points may be located on a great hyperbola branch as before;
but they can also be on a great ellipse, where it would be more reasonable to use trigonometric functions for defining a metric.
We have found the following definition to be the most logical:

\dfn{For every two points $x,y\in S^{1,1}$ we call
\[\dd(x,y)=\begin{cases}
 \arcosh(\llangle x,y\rrangle),&\mbox{if $x-y$ is timelike,}\\
 0,&\mbox{if $x-y$ is lightlike,}\\
 \infty,&\mbox{if } \llangle x,y\rrangle\leq-1 \mbox{ and } x\neq-y,\\
 \arccos(\llangle x,y\rrangle)&\mbox{otherwise}
\end{cases}\]
the \emph{(proper) de Sitter distance} of $x$ and $y$.}

This distance is well-defined, but obviously no metric.
Nevertheless, it is non-negative and symmetric,
and its value equals zero if and only if $x-y$ is located on the light cone.
Under certain circumstances, also the triangle inequality holds true;
but these circumstances will be dealt with in section \ref{Triangles}.

The cases in the definition also correspond to geometrical circumstances:
The first case ($x-y$ being timelike) is equivalent to $x$ and $y$ lying on a great hyperbola branch.
Vector $x-y$ lying on the light cone means, that the geodesic connecting $x$ and $y$ is a straight line.
In these cases, the de Sitter distance equals the ``time separation'' $\tau(x,y)$ (cf. O'Neill \cite{One}),
where $y$ has greater or equal $e_1$-coordinate.
In the remaining cases, this time seperation is zero.
The third case describes algebraically that $x$ and $y$ cannot be connected by a geodesic --
the reason is, that the two-dimensional subspace spanned by these vectors is timelike or lightlike
and $x$ and $y$ are located on the different components of the intersection with $S^{1,1}$.
And finally, what is left: The ``otherwise'' condition reflects the property of $x$ and $y$ to be points on a great ellipse.
Since the restriction of the Minkowski product to a spacelike vector subspace is a Riemannian metric,
the last case describes the distance function defined by this metric.

Finally, we subsume all distance functions under the term \emph{(generalized) de Sitter distance,}
symbolized with $d_\SS$ and defined as follows:
\[d_\SS(x,y)=\begin{cases}
 d_H(x,y),&\mbox{if } x,y\in H^2,\\
 d_H'(x,y),&\mbox{if } x,y\in (-H^2),\\
 \dd(x,y),&\mbox{if } x,y\in S^{1,1},\\
 \infty,&\mbox{otherwise.}
\end{cases}\]
This new distance function defines the distance for every pair of points in $\SS^2$.
The context will make clear, if the term ``de Sitter distance'' means the generalized or the proper de Sitter distance.

\vspace{.5cm}
For two different points $A,B\in\SS^2$ with $d_\SS(A,B)<\infty$ we find a tangent vector $X_{AB}$ at $A$ pointing in the direction of $B$ by computing
 \[X_{AB}:=\begin{cases} B-A,&\mbox{if $\spann\{A,B\}$ is lightlike,}\\
       \frac{B+\llangle A,B\rrangle A}{\llVert A\times B\rrVert},&\mbox{if $A,B\in(\pm H^2)$,}\\
       \frac{B-\llangle A,B\rrangle A}{\llVert A\times B\rrVert},&\mbox{otherwise.}
      \end{cases}\]

\Bem{Except for the first case, $X_{AB}$ is normalized.
 If $\spann(A,B)$ is spacelike or timelike, the denominator in the definition fraction can be replaced by $\sin(\dd(A,B))$ or $\sinh(d_\SS(A,B))$, respectively.}

Once we have tangent vectors, we can easily describe segments of de Sitter geodesics between two points.
Of course, if no such geodesic exists, the segment should be empty.
Thus, we define the \emph{(generalized) de Sitter segment} for $A\neq\pm B$, $d_\SS(A,B)<\infty$ as follows:
\[\ooverline{AB}:=\begin{cases} \{A+tX_{AB}\mid t\in[0,1]\},&\mbox{if $\spann\{A,B\}$ is lightlike,}\\
                   \{\cos(t)A+\sin(t)X_{AB}\mid t\in[0,\dd(A,B)]\},&\mbox{if $\spann\{A,B\}$ is spacelike,}\\
                   \{\cosh(t)A+\sinh(t)X_{AB}\mid t\in[0,d_\SS(A,B)]\},&\mbox{otherwise.}
                  \end{cases}\]

In analogy to the naming of de Sitter geodesics, we can distinguish hyperbolic, antipodal hyperbolic, and proper de Sitter segments.
In addition, proper de Sitter segments can be either line, great ellipse, or great hyperbola segments.
To not confuse great hyperbola segments on the de Sitter surface with hyperbolic segments,
We name proper de Sitter segments after the type of the plane whose intersection with $S^{1,1}$ gave the geodesic,
i.e.\ line segments are called lightlike, great ellipse segments spacelike, and great hyperbola segments timelike.

We have two remaining cases to consider:
If $A=B$, we define the generalized de Sitter segment to be
\[\ooverline{AB}:=\{A\};\]
and if $A=-B$ or $d_\SS(A,B)=\infty$, the segment is empty.

Note that for any Lorentz transformation $\Phi$,
\[\ooverline{\Phi(A)\Phi(B)}=\Phi\left(\ooverline{AB}\right)\]
holds.

\vspace{.5cm}
Finally, because we are going to obtain trigonometric laws, we still need to define angles.
We do this as follows:

\dfn{Let $A,B,C\in\SS^2$ be distinct.
 If $\ooverline{AB}$ and $\ooverline{AC}$ are of the same type (either hyperbolic, antipodal hyperbolic, proper de Sitter timelike, or spacelike segments) and neither lightlike nor empty,
 the de Sitter angle between these segments computes to
 \[\varangle(B,A,C):=\begin{cases}\arcosh(\lvert\llangle X_{AB},X_{AC}\rrangle\rvert),&\mbox{for proper de Sitter segments;}\\
                      \arccos(\llangle X_{AB},X_{AC}\rrangle),&\mbox{otherwise.}\end{cases}\]}

\Bem{Applying a Lorentz transformation does not change de Sitter angles.}

\Bem{We do not have to explicitly compute the relevant tangent vectors, because we find the relation
 \[\llangle X_{AB},X_{AC}\rrangle=\pm\leftllangle\frac{A\times B}{\llVert A\times B\rrVert},\frac{A\times C}{\llVert A\times C\rrVert}\rightrrangle,\]
 where the minus sign applies for $A,B,C\in S^{1,1}$, while the plus sign applies for $A,B,C\in H^2\cup(-H^2)$.}

\Beweis{Consider the case, that the relevant segments are proper de Sitter segments.
 The other case can be proved analogous.
 We have the definition
 \[X_{AB}=\frac{B-\llangle A,B\rrangle\cdot A}{\llVert A\times B\rrVert}.\]
 Solving the equation for $B$ and plugging the solution into $A\times B$ leads to
 \[A\times B=\llVert A\times B\rrVert\cdot A\times X_{AB}.\]
 Now, doing the same for $A\times C$ results in
 \begin{eqnarray*}
  \leftllangle \frac{A\times B}{\llVert A\times B\rrVert},\frac{A\times C}{\llVert A\times C\rrVert}\rightrrangle
  &=&\llangle A\times X_{AB},A\times X_{AC}\rrangle\\
  &=&-(\llangle A,A\rrangle\llangle X_{AB},X_{AC}\rrangle-\llangle A,X_{AC}\rrangle\llangle A,X_{AB}\rrangle)\\
  &=&-\llangle X_{AB},X_{AC}\rrangle.
 \end{eqnarray*}

 \vspace{-.3cm}}

\Bem{De Sitter angles are well-defined:
 If the involved segments are timelike, tangent vectors are timelike and normalized.
 WLOG, let $X_{AB}=e_1$. Denote the components of $X_{AC}$ by $x_1,x_2,x_3$. Then,
 \[|\llangle X_{AB},X_{AC}\rrangle|=|x_1|\geq\sqrt{x_1^2-x_2^2-x_3^2}=\llVert X_{AC}\rrVert=1\]
 holds.
 Thus, we can identify $|\llangle X_{AB},X_{AC}\rrangle|$ with a hyperbolic cosine.

 When dealing with (antipodal) hyperbolic segments, let WLOG be $A=\pm e_1$.
 Since the tangent vectors are Lorentz orthogonal to $A$, their first component must vanish, leading to
 \[\llangle X_{AB},X_{AC}\rrangle=\langle X_{AB},X_{AC}\rangle\leq\lVert X_{AB}\rVert \lVert X_{AC}\rVert=1,\]
 which allows to identify this product with a cosine.

 Finally, concerning spacelike segments, we note that $\JJ(A\times B)$ is Minkowski orthogonal to the (spacelike) vectors $A$ and $B$.
 Therefore it has to be timelike.
 Since $\JJ$ describes a Lorentz transformation, $A\times B$ is timelike as well.
 WLOG, we let $A\times B=e_1$ and continue in the same way as we did with timelike segments.}

\Bem{If $\ooverline{AB}$ and $\ooverline{AC}$ are of different types or both lightlike,
 the magnitude of the vertex angle at $A$ cannot reasonably be defined,
 for in this case the value of $\llangle X_B,X_C\rrangle$ can be any real number.}

\Bem{If one wants to define a congruence between angles on the de Sitter surface in a reasonable way,
 one must take into account the sign of $\llangle X_B,X_C\rrangle$.
 If one did not, there are angles obviously not being congruent, but having the same magnitude, e.g.
 $\varangle(B,A,B)$ and $\varangle(B,A,C)$ with $A$ being in $\ooverline{BC}\subset S^{1,1}$.}

\section{Classification of de Sitter Triangles} \label{Triangles}

By the term \emph{generalized de Sitter triangle} we denote any subset $\Delta\subset\SS^2$ consisting of three and only three elements.
This definition includes triangles that have empty sides.
It also includes degenerate triangles which have one vertex on the (open) segment between the other two,
but excludes the case of two vertices being identical.

\Konv{We name the vertices of a triangle using capital letters;
 the de Sitter distance or the segment between two vertices is denoted
 by the corresponding lower case variant of the labeling of the remainig point.}

We now want to classify generalized de Sitter triangles according to the types of their sides.

\dfn{A generalized de Sitter triangle $\Delta$ is said to be
 \begin{itemize}
  \item \DEF{hyperbolic,} if $\Delta\subset H^2$;
  \item \DEF{antipodal hyperbolic,} if $\Delta\subset(-H^2)$;
  \item \DEF{proper,} if $\Delta\subset S^{1,1}$; and
  \item \DEF{strange} otherwise.
 \end{itemize}
 The sides of strange triangles whose vertices are located on different components of $\SS^2$ are called 
 \DEF{strange segments.}
 Every de Sitter triangle that is not proper is \DEF{improper.}

 Proper de Sitter triangles may still have sides of different type.
 Therefore these can be further divided into the following:
 \begin{itemize}
  \item \DEF{spatiolateral triangles,} with all sides being spacelike;
  \item \DEF{chorosceles triangles,} with two sides being spacelike and the third being timelike;
  \item \DEF{tempolateral triangles,} with all sides being timelike;
  \item \DEF{chronosceles triangles,} with two sides being timelike and the third being spacelike;
  \item \DEF{lucilateral triangles,} with all sides being lightlike;
  \item \DEF{bimetrical chorosceles triangles,} with two sides being spacelike and the third being lightlike;
  \item \DEF{photosceles triangles with spacelike base,} with two sides being lightlike and the third being spacelike;
  \item \DEF{bimetrical chronosceles triangles,} with two sides being timelike and the third being lightlike;
  \item \DEF{photosceles triangles with timelike base,} with two sides being lightlike and the third being timelike; and
  \item \DEF{multiple triangles,} with one side of each type.
 \end{itemize}
 If at least one side of the triangle is empty, then $\Delta$ is called an \DEF{impossible triangle.}
 The empty sides are then called \DEF{impossible sides.}}

The names we gave to the different types of proper de Sitter triangles are meant to remind of equilateral and isosceles triangles, respectively.
The term ``bimetrical'' refers to the fact that the length of the lightlike segment is zero,
so only the two remaining sides of the triangle can be measured in a sensible way.

\vspace{.5cm}
When we start looking for examples of all these different types,
it comes out that lucilateral triangles are always degenerate.

To see this, we choose one vertex $A$ to equal $e_2$.
For lightlike segments, we have tangent vectors of the form $X_{AB}=B-A$.
These are Minkowski orthogonal to $A$ and lightlike, what leads to
$$X_{AB}=C_1\cdot(e_1\pm e_3),\quad X_{AC}=C_2\cdot(e_1\pm e_3),$$
where $C_{1,2}$ are arbitrary constants.
If $X_{AB}$ and $X_{AC}$ were not linearly dependent, we would easily obtain
that the side $a$ could only be lightlike, if one of the constants $C_{1,2}$ were zero, which cannot hold.
\hfill$\Box$

\begin{figure}
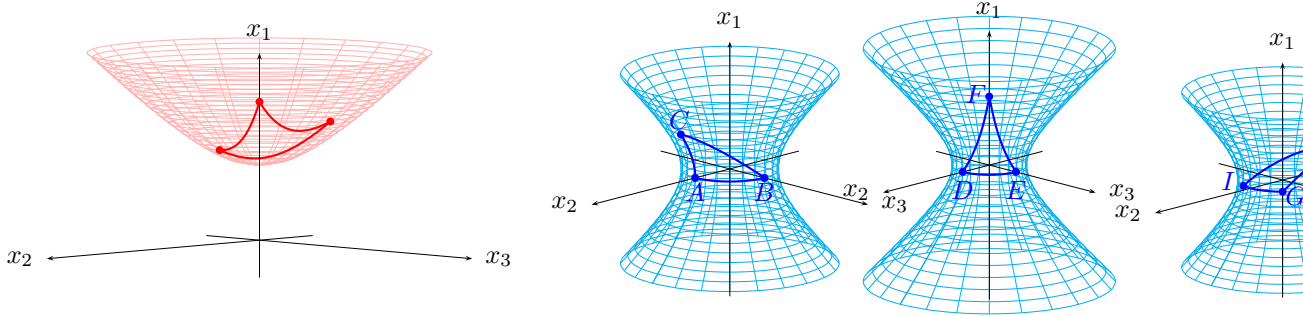

\begin{center}
   \pspicture(-3,-1)(3,1.25) 
    \psset{Beta=5,zMin=-0.5,zMax=2.5}
    \parametricplotThreeD[plotstyle=curve,hiddenLine=true,linecolor=lightred,linewidth=0.2pt](0,360)(1,2.5){%
     t cos u dup mul 1 neg add sqrt mul t sin u dup mul 1 neg add sqrt mul u}
    \parametricplotThreeD[plotstyle=curve,hiddenLine=true,linecolor=lightred,linewidth=0.2pt](1,2.5)(0,360){%
     u cos t dup mul 1 neg add sqrt mul u sin t dup mul 1 neg add sqrt mul t}

    \parametricplotThreeD[plotstyle=curve,linecolor=red](0,1.22513){0.75 t dup mul 1 add sqrt mul 0.03156 t mul add 0.99968 t mul 1.25 t dup mul 1 add sqrt mul 0.01895 t mul add}
    \parametricplotThreeD[plotstyle=curve,linecolor=red](0,1.37521){0.89058 t mul 1.33333 t dup mul 1 add sqrt mul 0.758 neg t mul add 1.66667 t dup mul 1 add sqrt mul 0.60641 neg t mul add}
    \parametricplotThreeD[plotstyle=curve,linecolor=red](0,1.82764){0.75 t dup mul 1 add sqrt mul 18.75 neg t mul 481 sqrt div add 16 t mul 481 sqrt div 1.25 t dup mul 1 add sqrt mul 11.25 neg t mul 481 sqrt div add}

    \pstThreeDCoor[nameX=$x_2$,nameY=$x_3$,nameZ=$x_1$,linewidth=0.2pt,linecolor=black]

    \pstThreeDDot[linecolor=red](1.22474,1.22474,2) 
    \pstThreeDDot[linecolor=red](0.75,0,1.25) 
    \pstThreeDDot[linecolor=red](0,1.33333,1.66667) 
   \endpspicture 
   \psset{unit=0.65cm}
   \pspicture(-5,-3)(3,1.25) 
    \psset{Beta=15,zMin=-2.7,zMax=2.7,xMin=-2,yMin=-2}
    \parametricplotThreeD[plotstyle=curve,hiddenLine=true,linecolor=cyan,linewidth=0.2pt](0,360)(-2,2){%
     t cos u dup mul 1 add sqrt mul t sin u dup mul 1 add sqrt mul u}
    \parametricplotThreeD[plotstyle=curve,hiddenLine=true,linecolor=cyan,linewidth=0.2pt](-2,2)(0,360){%
     u cos t dup mul 1 add sqrt mul u sin t dup mul 1 add sqrt mul t}

    \pstThreeDCoor[nameX=$x_2$,nameY=$x_3$,nameZ=$x_1$,linewidth=0.2pt,linecolor=black]

    \pstThreeDDot[linecolor=blue](1,0,0) \pstThreeDPut(1,0,-0.3){\blue$A$}
    \pstThreeDDot[linecolor=blue](0,1,0) \pstThreeDPut(0,1,-0.3){\blue$B$}
    \pstThreeDDot[linecolor=blue](1.41421,0,1) \pstThreeDPut(1.45,0,1.3){\blue$C$}
    \parametricplotThreeD[plotstyle=curve,linecolor=blue](0,90){t cos t sin 0}
    \parametricplotThreeD[plotstyle=curve,linecolor=blue](0,90){t cos 2 sqrt mul t sin t cos}
    \parametricplotThreeD[plotstyle=curve,linecolor=blue](0,1){t dup mul 1 add sqrt 0 t}
   \endpspicture 
   \psset{unit=0.5cm}
   \pspicture(-3,-4)(3,4.25) 
    \psset{Beta=15,zMin=-3.7,zMax=3.7,xMin=-2,yMin=-2}
    \parametricplotThreeD[plotstyle=curve,hiddenLine=true,linecolor=cyan,linewidth=0.2pt](0,360)(-3.2,3.2){%
     t cos u dup mul 1 add sqrt mul t sin u dup mul 1 add sqrt mul u}
    \parametricplotThreeD[plotstyle=curve,hiddenLine=true,linecolor=cyan,linewidth=0.2pt](-3.2,3.2)(0,360){%
     u cos t dup mul 1 add sqrt mul u sin t dup mul 1 add sqrt mul t}

    \pstThreeDCoor[nameX=$x_2$,nameY=$x_3$,nameZ=$x_1$,linewidth=0.2pt,linecolor=black]

    \pstThreeDDot[linecolor=blue](1,0,0) \pstThreeDPut(1,0,-0.5){\blue$D$}
    \pstThreeDDot[linecolor=blue](0,1,0) \pstThreeDPut(0,1,-0.5){\blue$E$}
    \pstThreeDDot[linecolor=blue](2,2,2.64575) \pstThreeDPut(2.5,2,2.8){\blue$F$}
    \parametricplotThreeD[plotstyle=curve,linecolor=blue](0,90){t cos t sin 0}
    \parametricplotThreeD[plotstyle=curve,linecolor=blue](0,1.73205){t dup mul 1 add sqrt 2 t mul 3 sqrt div 7 3 div sqrt t mul}
    \parametricplotThreeD[plotstyle=curve,linecolor=blue](0,1.73205){2 t mul 3 sqrt div t dup mul 1 add sqrt 7 3 div sqrt t mul}
   \endpspicture 
   \psset{unit=0.6cm}
   \pspicture(-4,-3)(3,3.25) 
    \psset{Beta=15,zMin=-2.7,zMax=2.7,xMin=-2,yMin=-2}
    \parametricplotThreeD[plotstyle=curve,hiddenLine=true,linecolor=cyan,linewidth=0.2pt](0,360)(-2,2){%
     t cos u dup mul 1 add sqrt mul t sin u dup mul 1 add sqrt mul u}
    \parametricplotThreeD[plotstyle=curve,hiddenLine=true,linecolor=cyan,linewidth=0.2pt](-2,2)(0,360){%
     u cos t dup mul 1 add sqrt mul u sin t dup mul 1 add sqrt mul t}

    \pstThreeDCoor[nameX=$x_2$,nameY=$x_3$,nameZ=$x_1$,linewidth=0.2pt,linecolor=black]

    {\psset{Alpha=90}
     \pstThreeDDot[linecolor=blue](1,0,0) \pstThreeDPut(1,0.3,-0.1){\blue$G$}
     \pstThreeDDot[linecolor=blue](1,1,1) \pstThreeDPut(1,1.3,1.1){\blue$H$}
     \pstThreeDDot[linecolor=blue](0.5,-0.86603,0) \pstThreeDPut(0.5,-1.2,0.1){\blue$I$}
     \parametricplotThreeD[plotstyle=curve,linecolor=blue](-60,0){t cos t sin 0}
     \parametricplotThreeD[plotstyle=curve,linecolor=blue](0,111.4707){t cos 0.9306 t sin mul add t cos 0.53729 t sin mul neg add t cos 0.39332 t sin mul add}
     \parametricplotThreeD[plotstyle=curve,linecolor=blue](0,1){1 t t}
    }
   \endpspicture 
   \psset{unit=0.6cm}
   \pspicture(-3,-3.75)(3,2.5) 
    \psset{Beta=15,zMin=-2.7,zMax=2.7,xMin=-2,yMin=-2}
    \parametricplotThreeD[plotstyle=curve,hiddenLine=true,linecolor=cyan,linewidth=0.2pt](0,360)(-2,2){%
     t cos u dup mul 1 add sqrt mul t sin u dup mul 1 add sqrt mul u}
    \parametricplotThreeD[plotstyle=curve,hiddenLine=true,linecolor=cyan,linewidth=0.2pt](-2,2)(0,360){%
     u cos t dup mul 1 add sqrt mul u sin t dup mul 1 add sqrt mul t}

    \pstThreeDCoor[nameX=$x_2$,nameY=$x_3$,nameZ=$x_1$,linewidth=0.2pt,linecolor=black]

    \pstThreeDDot[linecolor=blue](0.70711,0.70711,0) \pstThreeDPut(0.5,0.9,-0.1){\blue$J$}
    \pstThreeDDot[linecolor=blue](1.20711,0.20711,0.70711) \pstThreeDPut(1.4,0,0.8){\blue$K$}
    \pstThreeDDot[linecolor=blue](0,1.41421,1) \pstThreeDPut(-0.2,1.6,1.1){\blue$L$}
    \parametricplotThreeD[plotstyle=curve,linecolor=blue](0,1){1 t neg add 0.70711 mul 1 t add 0.70711 mul t}
    \parametricplotThreeD[plotstyle=curve,linecolor=blue](0,1){t 0.5 mul 0.70711 add t 0.5 mul neg 0.70711 add t 0.70711 mul}
    \parametricplotThreeD[plotstyle=curve,linecolor=blue](-0.41421,1){1.32623 t dup mul neg 1 add sqrt mul 0.87114 t dup mul neg 1 add sqrt mul 1.41421 t mul add 1.23198 t dup mul neg 1 add sqrt mul t add}
   \endpspicture 
   \psset{unit=0.5cm}
   \pspicture(-3,-4)(3,5) 
    \psset{Beta=15,zMin=-3.7,zMax=3.7,xMin=-3,yMin=-2}
    \parametricplotThreeD[plotstyle=curve,hiddenLine=true,linecolor=cyan,linewidth=0.2pt](0,360)(-3.2,3.2){%
     t cos u dup mul 1 add sqrt mul t sin u dup mul 1 add sqrt mul u}
    \parametricplotThreeD[plotstyle=curve,hiddenLine=true,linecolor=cyan,linewidth=0.2pt](-3.2,3.2)(0,360){%
     u cos t dup mul 1 add sqrt mul u sin t dup mul 1 add sqrt mul t}

    {\psset{Alpha=52}
     \pstThreeDCoor[nameX=$x_2$,nameY=$x_3$,nameZ=$x_1$,linewidth=0.2pt,linecolor=black]

     \pstThreeDDot[linecolor=blue](1,0,0) \pstThreeDPut(1,-0.2,-0.5){\blue$M$}
     \pstThreeDDot[linecolor=blue](1,1,1) \pstThreeDPut(0.8,1.35,1){\blue$N$}
     \pstThreeDDot[linecolor=blue](2,2.44949,-3) \pstThreeDPut(1.8,2.8,-3){\blue$O$}
     \parametricplotThreeD[plotstyle=curve,linecolor=blue](0,1){1 t t}
     \parametricplotThreeD[plotstyle=curve,linecolor=blue](0,1.73205){t dup mul 1 add sqrt 2 sqrt t mul 3 sqrt t mul neg}
     \parametricplotThreeD[plotstyle=curve,linecolor=blue](0,7.38207){0.73821 neg t mul t dup mul 1 add sqrt add 0.67732 neg t mul t dup mul 1 add sqrt add 1.41552 neg t mul t dup mul 1 add sqrt add}
    }
   \endpspicture 
   \psset{unit=0.6cm}
   \pspicture(-4,-3)(2.5,3.25) 
    \psset{Beta=15,zMin=-2.7,zMax=2.7,xMin=-2,yMin=-2}
    \parametricplotThreeD[plotstyle=curve,hiddenLine=true,linecolor=cyan,linewidth=0.2pt](0,360)(-2,2){%
     t cos u dup mul 1 add sqrt mul t sin u dup mul 1 add sqrt mul u}
    \parametricplotThreeD[plotstyle=curve,hiddenLine=true,linecolor=cyan,linewidth=0.2pt](-2,2)(0,360){%
     u cos t dup mul 1 add sqrt mul u sin t dup mul 1 add sqrt mul t}

    \pstThreeDCoor[nameX=$x_2$,nameY=$x_3$,nameZ=$x_1$,linewidth=0.2pt,linecolor=black]

    {\psset{Alpha=90}
     \pstThreeDDot[linecolor=blue](1,0,0) \pstThreeDPut(1,-0.3,-0.1){\blue$P_1$}
     \pstThreeDDot[linecolor=blue](1,1,1) \pstThreeDPut(1,1.3,1.1){\blue$P_2$}
     \pstThreeDDot[linecolor=blue](1,1,-1) \pstThreeDPut(1,1.3,-0.9){\blue$P_3$}
     \parametricplotThreeD[plotstyle=curve,linecolor=blue](0,1){1 t t}
     \parametricplotThreeD[plotstyle=curve,linecolor=blue](0,1){1 t t neg}
     \parametricplotThreeD[plotstyle=curve,linecolor=blue](0,2.82843){t dup mul 1 add sqrt 0.70711 t mul neg add t dup mul 1 add sqrt 0.70711 t mul neg add t dup mul 1 add sqrt 1.41421 t mul neg add}
    }
   \endpspicture 
   \psset{unit=0.6cm}
   \pspicture(-4,-3)(3,3.25) 
    \psset{Beta=15,zMin=-2.7,zMax=2.7,xMin=-2,yMin=-2}
    \parametricplotThreeD[plotstyle=curve,hiddenLine=true,linecolor=cyan,linewidth=0.2pt](0,360)(-2.2,2.2){%
     t cos u dup mul 1 add sqrt mul t sin u dup mul 1 add sqrt mul u}
    \parametricplotThreeD[plotstyle=curve,hiddenLine=true,linecolor=cyan,linewidth=0.2pt](-2.2,2.2)(0,360){%
     u cos t dup mul 1 add sqrt mul u sin t dup mul 1 add sqrt mul t}

    {\psset{Alpha=55}
     \pstThreeDCoor[nameX=$x_2$,nameY=$x_3$,nameZ=$x_1$,linewidth=0.2pt,linecolor=black]

     \pstThreeDDot[linecolor=blue](1,0,0) \pstThreeDPut(1,0,-0.3){\blue$Q_1$}
     \pstThreeDDot[linecolor=blue](1,2,2) \pstThreeDPut(1,2.5,2.1){\blue$Q_2$}
     \pstThreeDDot[linecolor=blue](0.44721,0.89443,0) \pstThreeDPut(0.44721,0.89443,-0.3){\blue$Q_3$}
     \parametricplotThreeD[plotstyle=curve,linecolor=blue](0,2){1 t t}
     \parametricplotThreeD[plotstyle=curve,linecolor=blue](0,63.43495){t cos t sin 0}
     \parametricplotThreeD[plotstyle=curve,linecolor=blue](0,2){t dup mul 1 add sqrt 0.44721 mul t dup mul 1 add sqrt 0.89443 mul t}
     }
   \endpspicture 
   \psset{unit=0.6cm}
   \pspicture(-2,-3)(3,3.25) 
    \psset{Beta=15,zMin=-2.7,zMax=2.7,xMin=-2,yMin=-2}
    \parametricplotThreeD[plotstyle=curve,hiddenLine=true,linecolor=cyan,linewidth=0.2pt](0,360)(-2.2,2.2){%
     t cos u dup mul 1 add sqrt mul t sin u dup mul 1 add sqrt mul u}
    \parametricplotThreeD[plotstyle=curve,hiddenLine=true,linecolor=cyan,linewidth=0.2pt](-2.2,2.2)(0,360){%
     u cos t dup mul 1 add sqrt mul u sin t dup mul 1 add sqrt mul t}

     \pstThreeDCoor[nameX=$x_2$,nameY=$x_3$,nameZ=$x_1$,linewidth=0.2pt,linecolor=black]

     \pstThreeDDot[linecolor=blue](1,0,0) \pstThreeDPut(0.8,-0.8,0){\blue$R$}
     \pstThreeDDot[linecolor=blue](2,1,2) \pstThreeDPut(2,1.5,2.1){\blue$S$}
     \pstThreeDDot[linecolor=blue](2,1,-2) \pstThreeDPut(2,1.5,-1.9){\blue$T$}
     \parametricplotThreeD[plotstyle=curve,linecolor=blue](0,1.73205){t dup mul 1 add sqrt 0.57735 t mul 1.1547 t mul}
     \parametricplotThreeD[plotstyle=curve,linecolor=blue](0,1.73205){t dup mul 1 add sqrt 0.57735 t mul 1.1547 t mul neg}
     \parametricplotThreeD[plotstyle=curve,linecolor=blue](0,8.94427){t dup mul 1 add sqrt 2 mul 1.78885 t mul neg add t dup mul 1 add sqrt 0.89443 t mul neg add t dup mul 1 add sqrt 2 mul 2.23607 t mul neg add}
   \endpspicture 
   \psset{unit=0.8cm}
   \pspicture(-2.7,-3)(2.7,3) 
    \psset{Beta=15,zMin=-2.7,zMax=2.7,xMin=-2,yMin=-2}
    \parametricplotThreeD[plotstyle=curve,hiddenLine=true,linecolor=cyan,linewidth=0.2pt](0,360)(-2.2,2.2){%
     t cos u dup mul 1 add sqrt mul t sin u dup mul 1 add sqrt mul u}
    \parametricplotThreeD[plotstyle=curve,hiddenLine=true,linecolor=cyan,linewidth=0.2pt](-2.2,2.2)(0,360){%
     u cos t dup mul 1 add sqrt mul u sin t dup mul 1 add sqrt mul t}

     \pstThreeDCoor[nameX=$x_2$,nameY=$x_3$,nameZ=$x_1$,linewidth=0.2pt,linecolor=black]

     \pstThreeDDot[linecolor=blue](1,0,0) \pstThreeDPut(1,0,-0.3){\blue$U$}
     \pstThreeDDot[linecolor=blue](0,1,0) \pstThreeDPut(0,1,-0.3){\blue$V$}
     \pstThreeDDot[linecolor=blue](0.71429,0.71429,0.14286) \pstThreeDPut(0.7,0.7,0.4){\blue$W$}
     \parametricplotThreeD[plotstyle=curve,linecolor=blue](0,90){t cos t sin 0}
     \parametricplotThreeD[plotstyle=curve,linecolor=blue](0,0.69985){t dup mul neg 1 add sqrt 1.02062 t mul 0.20412 t mul}
     \parametricplotThreeD[plotstyle=curve,linecolor=blue](0,0.69985){1.02062 t mul t dup mul neg 1 add sqrt 0.20412 t mul}
   \endpspicture 
   \psset{unit=0.6cm}
   \pspicture(-4,-3)(3,2.25) 
    \psset{Beta=15,zMin=-2.7,zMax=2.7,xMin=-2,yMin=-2}
    \parametricplotThreeD[plotstyle=curve,hiddenLine=true,linecolor=cyan,linewidth=0.2pt](135,315)(-2.2,2.2){%
     t cos u dup mul 1 add sqrt mul t sin u dup mul 1 add sqrt mul u}
    \parametricplotThreeD[plotstyle=curve,hiddenLine=true,linecolor=cyan,linewidth=0.2pt](-2.2,2.2)(135,315){%
     u cos t dup mul 1 add sqrt mul u sin t dup mul 1 add sqrt mul t}

    \pstThreeDCoor[nameX=$x_2$,nameY=$x_3$,nameZ=$x_1$,linewidth=0.2pt,linecolor=black]

    \pstThreeDDot[linecolor=blue](-0.71429,-0.71429,-0.14286) \pstThreeDPut(-0.7,-0.7,0.1){\blue$Z$}
    \parametricplotThreeD[plotstyle=curve,linecolor=blue](-0.71429,0.62923){t 1.02062 neg t dup mul neg 1 add sqrt mul 0.20412 neg t dup mul neg 1 add sqrt mul}
    \parametricplotThreeD[plotstyle=curve,linecolor=blue](-0.71429,0.62923){1.02062 neg t dup mul neg 1 add sqrt mul t 0.20412 neg t dup mul neg 1 add sqrt mul}

    \parametricplotThreeD[plotstyle=curve,hiddenLine=true,linecolor=cyan,linewidth=0.2pt](-45,135)(-2.2,2.2){%
     t cos u dup mul 1 add sqrt mul t sin u dup mul 1 add sqrt mul u}
    \parametricplotThreeD[plotstyle=curve,hiddenLine=true,linecolor=cyan,linewidth=0.2pt](-2.2,2.2)(-45,135){%
     u cos t dup mul 1 add sqrt mul u sin t dup mul 1 add sqrt mul t}

    \pstThreeDDot[linecolor=blue](1,0,0) \pstThreeDPut(1,0,-0.3){\blue$X$}
    \pstThreeDDot[linecolor=blue](0,1,0) \pstThreeDPut(0,1,-0.3){\blue$Y$}
    \parametricplotThreeD[plotstyle=curve,linecolor=blue](0,90){t cos t sin 0}
    \parametricplotThreeD[plotstyle=curve,linecolor=blue](0.62923,1){t 1.02062 neg t dup mul neg 1 add sqrt mul 0.20412 neg t dup mul neg 1 add sqrt mul}
    \parametricplotThreeD[plotstyle=curve,linecolor=blue](0.62923,1){1.02062 neg t dup mul neg 1 add sqrt mul t 0.20412 neg t dup mul neg 1 add sqrt mul}

   \endpspicture

\caption{Different types of de Sitter triangles: hyperbolic, chorosceles (1st row),
 chronosceles, bimetrical chorosceles, photosceles with spacelike base (2nd row),
 bimetrical chronosceles, photosceles with timelike base, multiple (3rd row),
 tempolateral, contractible and non-contractible spatiolateral (4th row) triangles} \label{Abb:Typen}
\end{center}
\end{figure}

For all of the other types, one can find non-degenerate examples (cf.~fig.~\ref{Abb:Typen}).

\vspace{.5cm}
Note that the last two examples in fig.~\ref{Abb:Typen} that have $U=X:=e_2$, $V=Y:=e_3$, and $W=-Z:=({\frac{1}{7}},{\frac{5}{7}},{\frac{5}{7}})^t$ are both spatiolateral,
but don't quite appear to be of same type.
It is useful here to make a distinction between these two types of spatiolateral triangles again:
\emph{contractible} triangles on the one hand, and \emph{non-contractible} ones on the other hand.

One can distinguish them by projecting onto $e_2$-$e_3$-plane,
where for contractible triangles it is possible to find a line through the origin
such that all three points lie on one side of the line,
whereas for non-contractible triangles this is impossible.
Another way of distinction is to compute the lengths of the sides of the triangle.
If they sum up to a number less than $2\pi$, the triangle is contractible,
and otherwise it is not. This will be proved in theorem \ref{Satz:gross_klein}.

\vspace{.5cm}
Let us finally investigate for what types of de Sitter triangles $\{A,B,C\}$ the triangle inequality
\[d_\SS(A,B)+d_\SS(B,C)\geq d_\SS(A,C)\]
holds.
Obviously, this is trivially true for strange triangles, because all of them have at least two strange sides.
The inequality also holds for (antipodal) hyperbolic triangles,
since $d_H$ is the standard metric and $d_H'$ is derived from it.
Photosceles triangles never satisfy the inequality (choose $\ooverline{AC}$ to be the base).
Bimetrical and multiple triangles satisfy the inequality if and only if they are isosceles
(which, in the case of multiple triangles, is mere chance).
For lucilateral triangles, the inequality trivially holds.
The same is true for non-contractible spatiolateral triangles, where the inequality follows from the fact
that the lengths of the sides sum up to a number greater than $2\pi$,
whereas each single side can only have a length less than $\pi$.
Impossible triangles that have two sides of infinite length satisfy the triangle inequality, too.
For impossible triangles that have only one such side, the inequality obviously does not hold.

There are even impossible triangles having no infinite side (namely those with two antipodal points).
For those, we find the triangle inequality to be true.
To see that, we start with noticing that they are just one case of degenerate triangles with all their vertices lying on a great ellipse.
Parametrizing this great ellipse allows us to write the distance of two vertices simply as a difference of two values for the parameter.
Having done this, we easily see the validity of the inequation.

The same method leads to the result, that every degenerate triangle satisfies the triangle inequality.

\satz{If the triangle $\Delta$ is chorosceles or chronosceles,
 then the triangle inequality generally does not hold.}

\Beweisenum{We give examples both for triangles that do not satisfy and for those that satisfy the inequality:
 \begin{enumerate}
  \item $\{A=e_2,B=e_3,C=({\frac{30}{41}\sqrt{2},\frac{59}{82}\sqrt{2},\frac{59}{82}\sqrt{2}})^t\}$
   is a chronosceles triangle. The lengths compute to
   \[a=b=\arcosh\left(\frac{59}{82}\sqrt{2}\right)\approx0.187;\,c=\frac{\pi}{2}\approx1.571,\]
   which shows that the length of the spacelike side $c$ is greater than the sum of the other lengths.
  \item For the chronosceles triangle $\{A,B,C\}$ where
   \[A=e_2,\quad B=\left({0,\frac{1}{2}\sqrt{3},-\frac{1}{2}}\right)^t,\quad C=({41,29,29})^t,\]
   we see that one of the timelike legs ($b$) is of greater length than the other two sides together.
  \item There are, however, some chronosceles triangles that satisfy the inequality,
   e.g.~$\{A=e_2,\, B=({0,\frac{1}{2}\sqrt{2},\frac{1}{2}\sqrt{2}})^t,\, C=({41,29,29})^t\}$.
  \item Now consider the chorosceles triangle $\{A,B,C\}$ given by
   \[A=({1,0,\sqrt{2}})^t,\quad B=({-1,0,\sqrt{2}})^t,\quad C=\left({0,\frac{1}{2}\sqrt{3},\frac{1}{2}}\right)^t.\]
   Its lengths compute to $a=b=\frac{\pi}{4}\approx0.785$ and $c=\arcosh(3)\approx1.763$,
   which gives $a+b<c$.
   It is the length of the timelike base that is too great.
  \item Another chorosceles triangle, one of whose legs is to great, is the following:
   \[\Delta=\left\{A=({1,0,\sqrt{2}})^t,\, B=\left({\frac{20}{21},0,\frac{29}{21}}\right)^t,\, C=\left({0,\cos\left(\frac{6\pi}{25}\right),\sin\left(\frac{6\pi}{25}\right)}\right)\right\}.\]
   Here, $a>b+c$ holds.
  \item Finally, we see that $\{A,B,C\}$ given by
   \[A=({1,0,\sqrt{2}})^t,\quad B=e_3,\quad C=\left({0,\frac{1}{2}\sqrt{3},\frac{1}{2}}\right)^t\]
   satisfies the triangle inequality. \hfill$\Box$
 \end{enumerate}}

\Lemma{\label{lemma:pos_neg_2}\label{bem:pos_neg_3}Every tempolateral triangle has one and only one vertex,
 at which the tangent vectors pointing in the directions of the other two vertices
 have different signs in their respective first component.}

\Beweis{Let $\Delta=\{A,B,C\}$ be tempolateral.
 To show the existence of such a vertex, we assume that the first component of $X_{AB}$ is positive.
 Then, by parametrizing the great hyperbola $A$ and $B$ are located on,
 we obtain the result that the first component of $X_{BA}$ must be negative.
 If $X_{BC}$ has a positive first component, the proof is done.
 Thus, let us assume that the first component of $X_{BC}$ is negative.
 Then, we find $X_{CB}$ having positive first component.
 Again, a negative first component of $X_{CA}$ would finish the proof,
 so we assume that $X_{CA}$ has a positive first component.
 But now we find $X_{AC}$ having a negative first component,
 giving the result that $A$ is a vertex with the desired property.

 Analogously one can see that assuming two such vertices implies
 that the tangent vectors at the third vertex must also have different signs in their respective first component.
 We obtain, thus, that the first component of $X_{AB}$ has the same sign as those of $X_{BC}$ and $X_{CA}$.
 Furthermore, we can assume that the first component of $A$ has the same sign, too.
 Then we compute
 \begin{eqnarray*}
  A&=&\cosh(b)C+\sinh(b)X_{{CA}}\\
  &=&\cosh(b)(\cosh(a)B+\sinh(a)X_{{BC}})+\sinh(b)X_{{CA}}\\
  &=&\cosh(b)(\cosh(a)(\cosh(c)A+\sinh(c)X_{{AB}})+\sinh(a)X_{{BC}})+\sinh(b)X_{{CA}}\\
  &=&\underbrace{\cosh(a)\cosh(b)\cosh(c)}_{>1}\cdot{}A+\underbrace{\cosh(a)\cosh(b)\sinh(c)}_{>0}\cdot{}X_{{AB}}\\
  &&\qquad+\underbrace{\sinh(a)\cosh(b)}_{>0}\cdot{}X_{{BC}}+\underbrace{\sinh(b)}_{>0}\cdot{}X_{{CA}},
 \end{eqnarray*}
 noticing that the first component of $A$ has increased (or decreased, depending on its sign),
 what cannot happen.
 Thus, we have shown the uniqueness of the vertex with the desired property.}

\Satz{\label{satz:Dreiecksungl_zeit} Non-degenerate tempolateral triangles do not satisfy the triangle inequality.}

\Beweis{Let $\{A,B,C\}$ be such a triangle with the first component of $X_{AB}$ and $X_{AC}$ having different signs.
 If we apply a Lorentz transformation that maps $X_{AB}$ to $e_1$, we easily see that
 \[\llangle X_{AB},X_{AC}\rrangle\geq1,\]
 where the equality holds for linearly dependent $X_{AB}$ and $X_{AC}$;
 but this would result in a degenerate triangle,
 so we obtain
 \begin{eqnarray*}
 \llangle X_{AB},X_{AC}\rrangle&>&1\\
 \sinh(b)\sinh(c)\llangle X_{AB},X_{AC}\rrangle&>&\sinh(b)\sinh(c).
 \end{eqnarray*}
 Finally, we compute
 \begin{eqnarray*}
  \cosh(a)&=&\llangle B,C\rrangle\\
  &=&\llangle\cosh(c)A+\sinh(c)X_{AB},\cosh(b)A+\sinh(b)X_{AC}\rrangle\\
  &=&\cosh(b)\cosh(c)+\sinh(b)\sinh(c)\llangle X_{AB},X_{AC}\rrangle\\
  &>&\cosh(b)\cosh(c)+\sinh(b)\sinh(c)\\
  &=&\cosh(b+c),
 \end{eqnarray*}
 which shows the relation
 \[a>b+c,\]
 contradicting the triangle inequality.}

\vspace{.5cm}
 \rm What is still left to show is the validity or invalidity of the triangle inequality for contractible spatiolateral triangles.
 We will show that the inequality does not hold in this case.
 This result, however, can most easily be obtained by using the duality between hyperbolic and de Sitter geometry,
 having us to postpone the proof to corollary \ref{Folg:gross_klein}.
 Actually, the proof for the validity of the triangle inequality for non-contractible spatiolateral triangles
 depends on the statement that the sum of lengths is greater than $2\pi$.
 Although that is a rather intuitive statement, it is still waiting for its proof,
 which will directly precede the corollary just mentioned.

\section{Polar Triangles}

The purpose of this section is to find a connection between hyperbolic triangles and proper de Sitter triangles.
Since we primarily want to investigate hyperbolic trigonometry,
we have the aim that triangles in $H^2\cup(-H^2)$ possess polar triangles in $S^{1,1}$,
even if not all the types of proper de Sitter triangles may occur.
The obvious way to provide such a mapping would be to consider the planes defining the sides of hyperbolic triangles,
then taking their intersection with de Sitter surface to obtain proper de Sitter geodesics,
and finally defining the intersection points of these geodesics as vertices of the polar triangle.
Unfortunately, these geodesics won't intersect.

Thus, we simply use the definition of polar triangles in spherical geometry and apply it to our case.

\dfn{The \emph{polar triangle} $\Delta'=\{A',B',C'\}$ of a generalized de Sitter triangle $\Delta=\{A,B,C\}$
 is defined by
 \[A':=\epsilon\cdot\frac{B\times C}{\llVert B\times C\rrVert},\quad B':=\epsilon\cdot\frac{C\times A}{\llVert C\times A\rrVert},\quad C':=\epsilon\cdot\frac{A\times B}{\llVert A\times B\rrVert},\]
 where $\epsilon=\sign(\det(A,B,C))=\pm1$.}

The polar triangle does not exist for every de Sitter triangle,
because one of the cross products may happen to be lightlike or zero.

If, on the other hand, $\Delta'$ exists, then it is well-defined,
and even the vertices of $\Delta'$ are each well-defined,
regardless of the order the vertices of $\Delta$ are given to find the polar triangle.
(Making this sure is the primary intent of $\epsilon$.)

Since the cross product $B\times C$ is orthogonal to $\spann(B,C)$,
one may think of vertex $A'$ as being determined (except for the sign) by the side $a$.

\vspace{.5cm}
It might appear more reasonable to define the vertices of the polar triangle
by using Lorentz orthogonality instead of simple Euclidean orthogonality.
If we attempted to do so, we could define a new ``Minkowski'' polar triangle $\Delta_{\rm Mink}'$ by
$A_{\mathrm{Mink}}'=\JJ(A'),\,B_{\mathrm{Mink}}'=\JJ(B'),\,C_{\mathrm{Mink}}'=\JJ(C');$
but this would only affect the actual position of the polar triangle,
whereas the lengths of its sides, its angles, and its type will be preserved.
So we stick to our first definition, because it is easier to handle.

\Bem{This definition satisfies our claim that any non-degenerate generalized de Sitter triangle $\Delta\subset H^2\cup(-H^2)$ has a polar triangle $\Delta'\subset S^{1,1}$.}

\Beweis{We know that $\JJ(A')\pperp B$.
A vector Minkowski orthogonal to a timelike one must be spacelike, so that is the case for $\JJ(A')$.
Hence, also $A'$ is spacelike and, since it is normalized, lies on $S^{1,1}$.}

\Satz{\label{Satz:Polar_Polar}Let $\Delta=\{A,B,C\}$ be a non-degenerate generalized de Sitter triangle.
If the polar triangle $\Delta'=\{A',B',C'\}$ exists, then the polar triangle of $\Delta'$ also exists, and
\[(A')'=A,\,(B')'=B,\,(C')'=C\]
hold.}

\Beweis{To see that $\Delta'$ is non-degenerate, we compute
 \begin{eqnarray*}
  \sign(\det(A',B',C'))&=&\sign(\det(\epsilon B\times C,\epsilon C\times A,\epsilon A\times B))\\
  &=&\epsilon\sign(\det(B\times C,C\times A,A\times B)))\\
  &=&\epsilon\sign(\langle(B\times C)\times(C\times A),A\times B\rangle)\\
  &=&\epsilon\sign(\langle B\times C,A\rangle\langle C\times A,B\rangle-\langle B\times C,B\rangle\langle C\times A,a\rangle)\\
  &=&\epsilon\sign(\det(A,B,C)^2)\\
  &=&\epsilon\\
  &\neq&0.
 \end{eqnarray*}
From the construction of $A'$ und $B'$ follows $A',B'\perp C$ (Euclidean orthogonality).
Furthermore, of course, $A'\times B'\perp A',B'$ holds.
Therefore $A'\times B'$ and $C$ are linearly dependent.
Particularly, $A'\times B'$ is of the same type as $C$ (either spacelike or timelike) and can thus be normalized.
Doing this, we get $(C')'=\pm C=:\delta C$.

Now we have on the one hand the relation
 \[\llangle A'\times B',(C')'\rrangle=\leftllangle A'\times B',\epsilon\frac{A'\times B'}{\llVert A'\times B'\rrVert}\rightrrangle=\epsilon\cdot(\sigma\llVert A'\times B'\rrVert),\]
 where $\sigma=+1$, if $A'\times B'$ is spacelike (and hence $C$ is spacelike, too),
 and $\sigma=-1$, if $A'\times B'$ (and hence $C$) is timelike.
 On the other hand we have:
 \begin{eqnarray*}
  \llangle A'\times B',(C')'\rrangle&=&\leftllangle\frac{B\times C}{\llVert B\times C\rrVert}\times\frac{C\times A}{\llVert C\times A\rrVert},\delta C\rightrrangle\\
  &=&\frac{\delta}{\llVert B\times C\rrVert\llVert C\times A\rrVert}\llangle \langle A,B\times C\rangle\cdot C-\underbrace{\langle C,B\times C\rangle}_{=0}\cdot A,C\rrangle\\
  &=&\frac{\delta}{\llVert B\times C\rrVert\llVert C\times A\rrVert}\det(A,B,C)\llangle C,C\rrangle\\
  &=&\frac{\delta\epsilon\cdot(\sigma\llVert C\rrVert)}{\llVert B\times C\rrVert\llVert C\times A\rrVert}.
 \end{eqnarray*}
Comparing these equations gives
\[\delta=\frac{1}{\llVert C\rrVert}\llVert A'\times B'\rrVert\llVert B\times C\rrVert\llVert C\times A\rrVert>0,\]
thus leading to $\delta=1$ and hence
 \[(C')'=C,\]
which was to show.}
\vspace{0.5cm}

The remaining part of this section deals with the connection between the type of a triangle and the type of its polar triangle.

\Satz{\label{Satz:kein_Polar} Bimetrical, photosceles, multiple, and lucilateral de Sitter triangles
 do not have a polar triangle, nor do such strange triangles that have two opposite points,
 or such impossible triangles one of whose impossible sides is contained in a lightlike plane.}

\Beweis{It is rather obvious that triangles with two opposing points have no polar triangles,
 for the cross product of those points equals zero.

 All the other cases describe a triangle that has at least one side located on a lightlike plane.
 The intersection of $S^{1,1}$ with such a plane is a pair of parallel lines, which can be parametrized as
 \[\{u\in\R^3\mid u=\pm A+t\cdot X,\,t\in\R\},\]
 where $A$ is one of the vertices on the lightlike or impossible side, and $X$ is Lorentz orthogonal to $A$.
 (In case of a lightlike side, this is a tangent vector).
 For a certain $t_0\in\R$, we thus have for the other vertex
 \[B=\pm A+t_0\cdot X,\]
 and we compute
 \[A\times B=A\times(\pm A+t_0\cdot X)=t_0\cdot A\times X.\]
 If $t=0$ (that means $B=-A$) or $A\times X=0$ holds, then it is again obvious that the polar triangle cannot exist.
 If, on the other hand, $t\cdot A\times X\neq0$ holds, we still know
 that $\JJ(A\times X)$ is Lorentz orthogonal to $A$ and $X$.
 If $\JJ(A\times X)\not\in\spann(A,X)$, we would have found a Minkowski orthogonal basis $\{A,X,\JJ(A\times X)\}$ of $\R^3$, which contains the lightlike vector $X$;
 but this is forbidden by lemma \ref{folg:Basis_zeit_2raum}.
 So we have $\JJ(A\times X)\in\spann(A,X)$.
 Because $\JJ(A\times X)$ is Lorentz orthogonal to any vector in $\spann(A,X)$,
 it is in particular Lorentz orthogonal to itself, meaning lightlike.
 Thus, $A\times B$ is also lightlike and therefore cannot be normalized.}

\Lemma{\label{lemma:exist_Polar} The polar triangle of a generalized de Sitter triangle $\Delta$
 which does not match one of the categories of theorem \ref{Satz:kein_Polar} always exists.}

\Beweis{The polar triangle exists if and only if the vectors $A\times B$,
 $B\times C$, and $C\times A$ can be normalized, i.e.~do not belong to the light cone.
 All triangles that contain two opposite points are dealt with in theorem \ref{Satz:kein_Polar} --
 strange triangles of that kind were mentioned explicitly,
 while two opposite points in proper de Sitter triangles would define an impossible side,
 which belongs to a plane of any type, and by this also to a lightlike plane.
 Therefore, the case that $A\times B=0$ is excluded.
 Furthermore, none of the planes $\spann(A,B)$, $\spann(A,C)$, and $\spann(B,C)$ is lightlike,
 so what we have to show is that no lightlike vector is Euclidean orthogonal to a spacelike or timelike plane.

 Again, we use the property of $\JJ$ to be among the Lorentz transformations
 to see that $A\times B$ is lightlike if and only if $\JJ(A\times B)$ is lightlike,
 and the latter is Lorentz orthogonal to $\spann(A,B)$.
 So we restate our aim to show as follows:
 There is no lightlike vector Minkowski orthogonal to a spacelike or timelike plane.

 Firstly, we notice that a lightlike vector contained in a timelike plane cannot be orthogonal to the same,
 because this would mean that there was a basis of that plane consisting of a spacelike and a lightlike vector,
 and thus the plane would be lightlike.

 So a Minkowski orthogonal basis of the plane $\spann(A,B)$ can be extended by $\JJ(A\times B)$
 to form a Minkowski orthogonal basis of $\R^3$.
 By Lemma \ref{folg:Basis_zeit_2raum}, we know that $\JJ(A\times B)$ is timelike if $\spann(A,B)$ is spacelike
 and vice versa, but never is it lightlike.}

\Satz{The polar triangle of a degenerate de Sitter triangle that does not match one of the categories of theorem\ref{Satz:kein_Polar}
 consists of only the zero vector.
 (It has never been stated that the polar triangle, if existing, has to be a de Sitter triangle!)}

\Beweis{By lemma \ref{lemma:exist_Polar}, we know that the polar triangle exists.
 We then compute
 \[\epsilon=\sign(\det(A,B,C))=\sign(0)=0\]
 and get the ``vertices'' $A'=B'=C'=0$.}

\Satz{\label{Satz:Polar_uneigtl_3raum}Improper de Sitter triangles on $H^2\cup(-H^2)$ without a pair of opposite points
 have spatiolateral polar triangles.}

\Beweis{Lemma \ref{lemma:exist_Polar} makes sure that $\Delta'=\{A',B',C'\}$ exists.
 Theorem \ref{Satz:Polar_Polar} says that $\Delta$ then is the polar triangle of $\Delta'$.
 In the proof of lemma \ref{lemma:exist_Polar} we saw that $A=(A')'$ is only timelike
 if $B'$ and $C'$ span a spacelike plane.
 Thus, since all of the points $A$, $B$, and $C$ are timelike, $\Delta'$ possesses only spacelike sides.}

\Satz{\label{folg:Polar_3raum_uneigtl}Polar triangles of non-degenerate spatiolateral de Sitter triangles
 are subsets of $H^2\cup(-H^2)$.}

\Beweis{The proof is analogous to the previous one.
 The only difference lies in the needed direction of an equivalence:
 $A'$ is timelike whenever $B$ and $C$ span a spacelike plane.}

\Bem{\label{bem:Polar_gross_hyp}Polar triangles of non-degenerate non-contractible spatiolateral triangles are not strange.
 To obtain this, name $A,B,C\in\Delta$ in such a way
 that moving on side $c$ from vertex $A$ to vertex $B$,
 then further along side $a$ to vertex $C$, and back on side $b$ to vertex $A$,
 means surrounding the $e_1$ axis in positive rotational direction.
 Since two sides alone will not add up to a full turn, each of the three sides must itself be passed in positive direction.
 Therefore, we know that each of the cross products $A\times B$, $B\times C$, and $C\times A$
 points in the direction of $e_1$, and thus, when normalized, lies on $H^2$.
 When searching for the vertices of the polar triangle,
 each of these normalized cross products are multiplied by the same $\epsilon=\pm1$.
 Hence, either are all of them in $(-H^2)$, or they all stay in $H^2$.

 Non-degenerate contractible spatiolateral triangles, on the other hand, all have strange polar triangles.
 If we had the sides from $A$ to $B$, from $B$ to $C$ and from $C$ to $A$ all moving in the same rotational direction around the $e_1$ axis,
 the triangle would surround this axis and could therefore not be contractible.
 Following the argumentation given above about non-contractible triangles,
 we see that the polar triangle must have both vertices in $H^2$ and vertices in $(-H^2)$.}

\Kor{\label{folg:Polar_hyp_gross}Polar triangles of (antipodal) hyperbolic triangles are non-contractible,
 whereas those of strange triangles with vertices in $H^2\cup(-H^2)$ are contractible.}

\Satz{\label{folg:Polar_schenklig_seltsam}Chronosceles and chorosceles triangles both have strange polar triangles.
 The same is true for impossible triangles that contain a spacelike side and do not match one of the categories of theorem \ref{Satz:kein_Polar}.}

\Beweis{We already saw in the proof of lemma \ref{lemma:exist_Polar} that spacelike sides in a triangle
 correspond to timelike vertices in the polar triangle and vice versa.
 The statement about impossible triangles results from the fact
 that the vertices of the impossible sides span a timelike plane,
 or otherwise the triangle would be mentioned by theorem \ref{Satz:kein_Polar}.}

\Satz{\label{Satz:Polar_3zeit_unmoegl_uu} Any non-degenerate tempolateral triangle
 possesses an impossible polar triangle with one timelike side and two impossible sides,
 and each of these impossible sides is defined by a timelike plane.
 Concerning this kind of impossible triangles, the corresponding polar triangle is either of the same type
 or tempolateral.}

\Beweis{The polar triangle of a tempolateral triangle consists of only spacelike vertices.
 We compute
 \begin{eqnarray*}
  \llangle A',B'\rrangle&=&\leftllangle\epsilon\frac{B\times C}{\llVert B\times C\rrVert},\epsilon\frac{C\times A}{\llVert C\times A\rrVert}\rightrrangle\\
  &=&\underbrace{\epsilon^2}_{=1}\llangle X_{CA},X_{CB}\rrangle,
 \end{eqnarray*}
 and in analogy
 \[\llangle B',C'\rrangle=\llangle X_{AB},X_{AC}\rrangle,\quad\llangle A',C'\rrangle=\llangle X_{BA},X_{BC}\rrangle.\]
 Let $X_{AB},X_{AC}$ be the pair of tangent vectors at a vertex that have different signs in their first components.
 We already noticed in theorem \ref{satz:Dreiecksungl_zeit} that under this condition,
 $\llangle B',C'\rrangle=\llangle X_{AB},X_{AC}\rrangle>1$ holds.
 This means that $B'$ and $C'$ span a timelike plane.
 (One can easily parametrize a great ellipse or a pair of straight lines to see
 that two points on such a de Sitter geodesic have Minkowski product of norm less or equal than 1.)
 Furthermore, $B'$ and $C'$ are on the same branch of the great hyperbola,
 because otherwise the Minkowski product would be negative.
 Thus, $a'$ is timelike.
 By the same argumentation, we have for the other sides
 \[\llangle A',B'\rrangle<-1 \mbox{ and } \llangle A',C'\rrangle<-1,\]
 which describes two impossible sides.

 Now let $\{A,B,C\}$ be an impossible triangle as described above.
 Let $a$ and $b$ be the impossible sides.
 That means, $B$ and $C$ span a timelike plane, but lie on the different branches of the corresponding great hyperbola,
 and the same holds true for $C$ and $A$.
 Now, replace $C$ by $-C$.
 We have a new triangle $\{\oldA,\oldB,\oldC\}$ with $\oldA=A,\,\oldB=B,\,\oldC=-C,$ which is timelike.
 One can easily verify that $\oldA'=A'$, $\oldB'=B'$, and $\oldC'=-C'$ hold.
 If $\oldA$ is the vertex with different signs in the first components of the tangent vectors,
 we have $\olda'$ being timelike and $\oldb',\,\oldc'$ both being impossible.
 Now re-substitute $C'$ for $\oldC'=-C'$, whence $c'=\oldc'$ is not affected at all,
 $b'$ becomes timelike and $a'$ becomes impossible.
 So we got a polar triangle of the same type as the original triangle.
 No difference appears if we assume that $\oldB$ is the vertex with opposing tangent vectors.
 If $\oldC$ is this vertex, $\oldc'$ is timelike, whereas $\oldb'$ and $\olda'$ are impossible.
 Re-substituting $C$ for $\oldC'$, we do not notice any influence on $c'=\oldc'$,
 but both impossible sides $\olda',\oldb'$ change for timelike sides $a',b'$.
 Thus, the polar triangle $\{A',B',C'\}$ in this case is timelike.}

\Satz{\label{Satz:Polar_unmoegl_unmoegl} Let $\Delta$ be an impossible triangle that does not match
 one of the categories of theorem \ref{Satz:kein_Polar}.
 If $\Delta$ has two timelike sides, the same holds for $\Delta'$.
 If all the sides of $\Delta$ are impossible, the same is true for $\Delta'$.}

\Beweis{Let $b$ and $c$ be two lightlike sides of $\Delta$, and $a$ be impossible.
 If the tangent vectors $X_{AB}$ and $X_{AC}$ had different signs in their respective first component,
 then we know from the proof of the triangle inequality for lucilateral triangles (theorem \ref{satz:Dreiecksungl_zeit}), that
 \[\llangle B,C\rrangle>\cosh(b+c)>1\]
 holds.
 However, this inequality describes the property of $a$ being timelike and not impossible.
 Thus, we have the first components  $X_{AB}$ and $X_{AC}$ bearing the same sign,
 which gives us an impossible side $a'$, according to the previous proof.

 Replacing $C$ with $-C$ gives a new triangle $\{\oldA=A,\oldB=B,\oldC=-C\}$,
 that has two timelike sides $\olda$ and $\oldc$, and one impossible side $\oldb$.
 Following our recent thoughts, $\oldb'$ has to be impossible.
 Re-substituting $C'$ for $\oldC'=-C'$ gives a timelike side $b'$.
 By the same argumentation, we have $a'$ being timelike.

 Now, if all the sides of $\Delta$ are impossible, we replace $A$ with $-A$ to get the triangle
 $\{\oldA=-A,\oldB=B,\oldC=C\}$,
 which has two timelike sides $\oldb$ and $\oldc$, and still one impossible side $\olda$.
 As we already know, under this conditions $\olda'$ is impossible,
 whereas $\oldb'$ and $\oldc'$ are timelike.
 Re-substituting $A'$ by $\oldA'=-A'$ then gives three impossible sides $a',$ $b',$ and $c'$.}

\Satz{Finally, the polar triangle of a non-degenerate strange triangle that has a point in $S^{1,1}$
 is either of the same type, chorosceles, chronosceles, or impossible with one spacelike side.
 In the latter case, the polar triangle does not match one of the categories of theorem \ref{Satz:kein_Polar}.}

\Beweis{Firstly, consider the case that the triangle has a spacelike side $a$, and vertex $A\in H^2\cup(-H^2)$.
 Then we have $A'\in H^2\cup(-H^2)$ and $B',C'\in S^{1,1}$,
 which means we get a polar triangle that is strange, non-degenerate and no subset of $H^2\cup(-H^2)$.

 The other cases result from the theorems \ref{folg:Polar_schenklig_seltsam} and \ref{Satz:Polar_Polar}.}

\section{Trigonometry of $\SS^2$}

This section gives an application of polar triangles.
We investigate the relations between the angles and the lengths of sides in a triangle.
To this end, we have to be able to measure these quantities,
which is only possible if the sides are spacelike or timelike (to measure their lengths)
and if adjacent sides are of same type (to measure the angle between them).
For a degenerate triangle, angles do not contain much information,
so we restrict our analysis to non-degenerate triangles that are either spatiolateral, tempolateral, or (antipodal) hyperbolic.

We abbreviate the term \emph{law of cosines for sides} by LCS and \emph{law of cosines for angles} by LCA.

All proofs in this section are rather elementary,
i.e.\ they contain no ``tricks'' or unexpected transformations.

\subsection*{Laws of Cosines}

\satz[Hyperbolic and Antipodal Hyperbolic LCS]{%
 Let $\Delta=\{A,B,C\}$ be an (antipodal) hyperbolic triangle.
 Denote the sides as usual by $a$, $b$, and $c$,
 denote the angle at vertex $A$ by $\varangle(B,A,C)=:\alpha$, at vertex $B$ by $\beta$,
 and at vertex $C$ by $\gamma$.
 Then,
 \[\cosh(a)=\cosh(b)\cosh(c)-\cos(\alpha)\sinh(b)\sinh(c),\]
 \[\cosh(b)=\cosh(a)\cosh(c)-\cos(\beta)\sinh(a)\sinh(c), \makebox[0pt][l]{ and}\]
 \[\cosh(c)=\cosh(a)\cosh(b)-\cos(\gamma)\sinh(a)\sinh(b)\]
 hold.}

\rm This theorem can be found in most works about hyperbolic geometry (see the references section)
with proofs of different degrees of difficulty.
A proof similar to ours is given by Iversen \cite{Ive},
but he does the computation in a more sophisticated vector space (the $sl_2$ space).

\Beweis{We obtain this result by simply computing
\begin{eqnarray*}
 \cos(\alpha)&=&\llangle X_{AB},X_{AC}\rrangle\\
 &=&\leftllangle\frac{B-\cosh(c)A}{\sinh(c)},\frac{C-\cosh(b)A}{\sinh(b)}\rightrrangle,
\end{eqnarray*}
hence
\begin{eqnarray*}
 \cos(\alpha)\sinh(b)\sinh(c)&=&\llangle B,C\rrangle-\cosh(c)\cdot\llangle A,C\rrangle-\cosh(b)\cdot\llangle B,A\rrangle\\
  &&\qquad+\cosh(b)\cosh(c)\cdot\llangle A,A\rrangle\\
 &=&-\cosh(a)+\cosh(b)\cosh(c).
\end{eqnarray*}
Renaming the vertices yields the other equations.}

\satz[LCS in Non-Contractible Spatiolateral Triangles]{%
 For any non-contractible spatiolateral de Sitter triangle with sides and angles named as previously,
 the following equalities hold.
 \[\cos(a)=\cos(b)\cos(c)-\cosh(\alpha)\sin(b)\sin(c),\]
 \[\cos(b)=\cos(a)\cos(c)-\cosh(\beta)\sin(a)\sin(c),\]
 \[\cos(c)=\cos(a)\cos(b)-\cosh(\gamma)\sin(a)\sin(b).\]}

\Beweis{As we already saw in the proof of theorem \ref{Satz:Polar_3zeit_unmoegl_uu},
\[\llangle X_{AB},X_{AC}\rrangle=\llangle B',C'\rrangle\]
holds.
This expression is negative, because $B'$ and $C'$ are either both hyperbolic or both antipodal hyperbolic.
Thus, we obtain
 \[\cosh(\alpha)=\lvert\llangle X_{AB},X_{AC}\rrangle\rvert=-\llangle X_{AB},X_{AC}\rrangle,\]
and hence
 \[-\cosh(\alpha)\sin(b)\sin(c)=\cos(a)-\cos(b)\cos(c).\]

 \vspace{-.3cm}}

\satz[LCS in Contractible Spatiolateral Triangles]{%
 Let $\{A,B,C\}$ be a contractible spatiolateral de Sitter triangle.
 Name the vertices in such a way that $a'$ is not strange.
 Then, the following is true.
 \[\cos(a)=\cos(b)\cos(c)-\cosh(\alpha)\sin(b)\sin(c),\]
 \[\cos(b)=\cos(a)\cos(c)+\cosh(\beta)\sin(a)\sin(c),\]
 \[\cos(c)=\cos(a)\cos(b)+\cosh(\gamma)\sin(a)\sin(b).\]}

The proof is completely analogous to the previous proofs.

\satz[LCS in Tempolateral Triangles]{%
 Let $\Delta$ be a tempolateral de Sitter triangle
 with $X_{AB}$ and $X_{AC}$ having different signs in their respective first component.
 Then,
 \[\cosh(a)=\cosh(b)\cosh(c)+\cosh(\alpha)\sinh(b)\sinh(c),\]
 \[\cosh(b)=\cosh(a)\cosh(c)-\cosh(\beta)\sinh(a)\sinh(c), \makebox[0pt][l]{ and}\]
 \[\cosh(c)=\cosh(a)\cosh(b)-\cosh(\gamma)\sinh(a)\sinh(b)\]
 hold.}

\Beweis{We have proven in theorem \ref{satz:Dreiecksungl_zeit} that
\[\llangle X_{AB},X_{AC}\rrangle>1>0\]
holds.
For the other vertices, we have
\[\llangle X_{BC},X_{BA}\rrangle<0 \quad\mbox{and}\quad \llangle X_{CA},X_{CB}\rrangle<0\]
instead.
What remains is a computation that is again completely analogous to the previous ones.}

\Lemma{\label{Satz:H2_S11} Let $\{A,B,C\}$ be (antipodal) hyperbolic with sides and angles named as usual.
 Let $\{A',B',C'\}$ be the corresponding polar triangle with accordingly named sides and angles.
 Then we have the following correlation between the sides and angles of these triangles.
 \[\alpha=\pi-a',\;\beta=\pi-b',\;\gamma=\pi-c';\]
 \[a=\alpha',\;b=\beta',\;c=\gamma'.\]}

\Beweis{We already know that
\[\cos(\alpha)=\llangle X_{AB},X_{AC}\rrangle=-\llangle B',C'\rrangle=-\cos(a')\]
holds, since the polar triangle is spatiolateral according to theorem \ref{Satz:Polar_uneigtl_3raum}.
With $\alpha,a'\in[0,\pi]$ we obtain the desired result $a'=\pi-\alpha$, or $\alpha=\pi-a'$.
Furthermore, for the polar triangle, which is non-contractible, we have
\[-\cosh(\alpha')=\llangle X',Y'\rrangle=\llangle B,C\rrangle=-\cosh(a),\]
from which $\alpha'=a$ results.}

\satz[Hyperbolic and Antipodal Hyperbolic LCA]{\label{hyp_Winkelcos}
 Let $\{A,B,C\}$ be an (antipodal) hyperbolic triangle with sides and angles named as usual. Then,
 \[\cos(\alpha)=-\cos(\beta)\cos(\gamma)+\cosh(a)\sin(\beta)\sin(\gamma),\]
 \[\cos(\beta)=-\cos(\alpha)\cos(\gamma)+\cosh(b)\sin(\alpha)\sin(\gamma), \makebox[0pt][l]{ and}\]
 \[\cos(\gamma)=-\cos(\alpha)\cos(\beta)+\cosh(c)\sin(\alpha)\sin(\beta)\]
 hold.}

This theorem is also mentioned by Anderson \cite{And}, Iversen \cite{Ive}, and Thurston \cite{Thu},
but they all give considerable longer proofs.
Although Thurston uses a similar duality like we do, his proof stays rather complex
since the hyperbolic laws of cosines are derived by means of Minkowski product matrices.

\Beweis{The polar triangle $\{A',B',C'\}$ is non-contractible spatiolateral.
 By the LCS for triangles of this type, we have
 \[\cos(a')=\cos(b')\cos(c')-\cosh(\alpha')\sin(b')\sin(c').\]
 Lemma \ref{Satz:H2_S11} leads to the desired result.}

\satz[LCA in Non-Contractible Spatiolateral Triangles]{%
 For any non-contractible spatiolateral de Sitter triangle with the usual labels,
 \[\cosh(\alpha)=\cosh(\beta)\cosh(\gamma)+\cos(a)\sinh(\beta)\sinh(\gamma),\]
 \[\cosh(\beta)=\cosh(\alpha)\cosh(\gamma)+\cos(b)\sinh(\alpha)\sinh(\gamma), \makebox[0pt][l]{ and}\]
 \[\cosh(\gamma)=\cosh(\alpha)\cosh(\beta)+\cos(b)\sinh(\alpha)\sinh(\beta)\]
 hold.}

\Beweis{Now, the polar triangle is hyperbolic or antipodal hyperbolic.
 From the (antipodal) hyperbolic LCS we know that
 \[\cosh(a')=\cosh(b')\cosh(c')-\cos(\alpha')\sinh(b')\sinh(c')\]
 holds.
 Lemma \ref{Satz:H2_S11}, applied to $\Delta'$, leads to the equations above.}

\satz[LCA in Contractible Spatiolateral Triangles]{\label{klein_Winkelcos}
 Let $\{A,B,C\}$ be a contractible spatiolateral de Sitter triangle with side $a'$ of the polar triangle not being strange.
 Then the following equations hold.
 \[\cosh(\alpha)=\cosh(\beta)\cosh(\gamma)+\cos(a)\sinh(\beta)\sinh(\gamma),\]
 \[\cosh(\beta)=\cosh(\alpha)\cosh(\gamma)-\cos(b)\sinh(\alpha)\sinh(\gamma),\]
 \[\cosh(\gamma)=\cosh(\alpha)\cosh(\beta)-\cos(c)\sinh(\alpha)\sinh(\beta).\]}

\Beweis{We change the orientation of vertex $A'$ in the polar triangle
 to get the hyperbolic or antipodal hyperbolic triangle $\{\oldA'=-A',\oldB'=B',\oldC'=C'\}$.
 Now let us look for a relation between the angles of $\{A,B,C\}$ and the sides of $\{\oldA',\oldB',\oldC'\}$:
 \[-\cosh(\alpha)=\llangle X_{AB},X_{AC}\rrangle=\llangle B',C'\rrangle=\llangle\oldB',\oldC'\rrangle=-\cosh(\olda'),\]
 \[\cosh(\beta)=\llangle X_{BA},X_{BC}\rrangle=\llangle A',C'\rrangle=-\llangle\oldA',\oldC'\rrangle=\cosh(\oldb'),\]
 and in analogy $\cosh(\gamma)=\cosh(\oldc')$.
 That resembles the relations we already know from non-contractible triangles.
 But what about the sides of $\{A,B,C\}$ and the angles $\alpha',\beta',\gamma'$ of $\{\oldA',\oldB',\oldC'\}$?
 Since we know that $\{\oldA=-A,\oldB=B,\oldC=C\}$ is the polar triangle of $\{\oldA',\oldB',\oldC'\}$,
 we can compute
 \[\cos(a)=\llangle B,C\rrangle=\llangle\oldB,\oldC\rrangle=-\llangle X_{\oldA'\oldB'},X_{\oldA'\oldC'}\rrangle=-\cos(\alpha'),\]
 \[\cos(b)=\llangle A,C\rrangle=-\llangle\oldA,\oldC\rrangle=\llangle X_{\oldB'\oldA'},X_{\oldB'\oldC'}\rrangle=\cos(\beta'),\]
 and $\cos(c)=\cos(\gamma')$.
 The LCS for $\{\oldA',\oldB',\oldC'\}$ yields the desired result.}

\satz[LCA in Tempolateral Triangles]{\label{zeit_Winkelcos}
 Let $\Delta$ be tempolateral
 with the tangent vectors at $A$ bearing different signs in their respective first component.
 Then we have the following relations.
 \[\cosh(\alpha)=\cosh(\beta)\cosh(\gamma)+\cosh(a)\sinh(\beta)\sinh(\gamma),\]
 \[\cosh(\beta)=\cosh(\alpha)\cosh(\gamma)-\cosh(b)\sinh(\alpha)\sinh(\gamma),\]
 \[\cosh(\gamma)=\cosh(\alpha)\cosh(\beta)-\cosh(c)\sinh(\alpha)\sinh(\beta).\]}

\Beweis{We replace the vertex $A'$ of the polar triangle by $-A'$ and get the tempolateral triangle
 $\{\oldA'=-A',\oldB'=B',\oldC'=C'\}$ with different signs in the first components of the tangent vectors at $\oldA'$.
 Following the previous proof and applying the $LCS$ to $\{\oldA',\oldB',\oldC'\}$,
 we get the desired equations.}
\vspace{0.5cm}

\rm Let us close our investigation of the laws of cosines by considering the case of degenerate triangles.
How do our achievements appear now?
For the LCS's, we simply have rules of the type
\[\cos(a)=\cos(b+c).\]
Again, we have proven that the triangle inequality holds for degenerate triangles.

Regarding the LCA's, we notice that the sine terms vanish, whence in most cases we have vacant formulae like this:
\[\cos(0)=\cos(0)\cdot\cos(0),\]
optionally with hyperbolic cosine.
The only more ``interesting'' result appears from the (antipodal) hyperbolic LCA, which gives
\[\cos(\alpha)=-\cos(\beta)\cos(\gamma).\]
The cosines can only take the values $-1$ or $1$,
thus we have either all angles amount to $\pi$ or one angle amount to $\pi$ and the other ones vanish.
The first case can be eliminated in analogy to the second part of the proof of lemma \ref{lemma:pos_neg_2}.

To conclude, the only case with any mathematical content shows us the hardly astonishing fact
that a degenerate (antipodal) hyperbolic triangle has one straight and two zero angles.
By now, it should be obvious why we can constrict our investigation to non-degenerate triangles.

\subsection*{Laws of Sines} \enlargethispage{2\baselineskip}

With the concept of polar triangles, we quickly deduce the laws of sines.

\satz[Hyperbolic and Antipodal Hyperbolic Law of Sines]{%
 For any (antipodal) hyperbolic triangle with labels as usual,
 \[\frac{\sin(\alpha)}{\sinh(a)}=\frac{\sin(\beta)}{\sinh(b)}=\frac{\sin(\gamma)}{\sinh(c)}\]
 holds.}

\rm This theorem is subject of the most books about hyperbolic geometry
(of the literature mentioned in the references section, only Ungar \cite{Ung} does not mention it).
We found none of the proofs being similar to ours.

\Beweis{We know that
 \begin{eqnarray*}
  \sin(\gamma)&=&\sin(\pi-c')=\sin(c')\\
  &=&\llVert A'\times B'\rrVert\\
  &=&\leftllVert\epsilon\frac{B\times C}{\llVert B\times C\rrVert}\times\left(\epsilon\frac{C\times A}{\llVert C\times A\rrVert}\right)\rightrrVert\\
  &=&\frac{1}{\llVert B\times C\rrVert\llVert C\times A\rrVert}\llVert\langle B,C\times A\rangle\cdot C-\langle C,C\times A\rangle\cdot B\rrVert\\
  &=&\frac{\lvert\det(A,B,C)\rvert}{\llVert B\times C\rrVert\llVert C\times A\rrVert},
 \end{eqnarray*}
 and thus
 \[\frac{\sin(\gamma)}{\sinh(c)}=\frac{\lvert\det(A,B,C)\rvert}{\llVert A\times B\rrVert\llVert B\times C\rrVert\llVert C\times A\rrVert}.\]
 This expression is entirely symmetric in $A$, $B$, and $C$,
 whence it remains the same for the other fractions
 \[\frac{\sin(\alpha)}{\sinh(a)}\mbox{ and }\frac{\sin(\beta)}{\sinh(b)}.\]

 \vspace{-.3cm}}

\satz[Law of Sine in Non-Contractible Spatiolateral Triangles]{%
 For any non-contractible spatiolateral triangle labeled as usual,
 \[\frac{\sinh(\alpha)}{\sin(a)}=\frac{\sinh(\beta)}{\sin(b)}=\frac{\sinh(\gamma)}{\sin(c)}\]
 holds.}

\Beweis{This follows from $\{A',B',C'\}$ being (antipodal) hyperbolic and
 \[\frac{\sinh(\alpha)}{\sin(a)}=\frac{\sinh(a')}{\sin(\alpha')}.\]

 \vspace{-.3cm}}

\satz[Law of Sine in Contractible Spatiolateral Triangles]{%
 Let $\{A,B,C\}$ be such a triangle with the side $a'$ of the polar triangle not being strange.
 Then the following equation holds.
 \[-\frac{\sinh(\alpha)}{\sin(a)}=\frac{\sinh(\beta)}{\sin(b)}=\frac{\sinh(\gamma)}{\sin(c)}.\]}

\Beweis{By changing the orientation of $A'$ we get the (antipodal) hyperbolic triangle $\{\oldA',\oldB',\oldC'\}$
 (cf.\ theorem \ref{klein_Winkelcos}).
 Now we compute
 \begin{eqnarray*}
  -\frac{\sinh(\alpha)}{\sin(a)}=-\frac{\sinh(\olda')}{-\sin(\alpha')}&=&\frac{\sinh(\oldb')}{\sin(\beta')}=\frac{\sinh(\beta)}{\sin(b)}\\
  &=&\frac{\sinh(\oldc')}{\sin(\gamma')}=\frac{\sinh(\gamma)}{\sin(c)}.
 \end{eqnarray*}

 \vspace{-.3cm}}

\satz[Law of Sines in Tempolateral Triangles]{%
 Let $\Delta$ be a tempolateral triangle labeled as usual.
 Let the first components of the tangent vectors at $A$ have different signs.
 Then,
 \[-\frac{\sinh(\alpha)}{\sinh(a)}=\frac{\sinh(\beta)}{\sinh(b)}=\frac{\sinh(\gamma)}{\sinh(c)}\]
 holds.}

\Beweis{In the proof of theorem \ref{zeit_Winkelcos} we constructed various auxiliary triangles,
 to which we gave blackletter labels.
 These triangles are still useful for this proof, when we compute
 \begin{eqnarray*}
  \sinh(\gamma)&=&\sinh(\oldc')\\
  &=&\llVert\oldA'\times\oldB'\rrVert\\
  &=&-\llVert A'\times B'\rrVert\\
  &=&-\frac{\lvert\det(A,B,C)\rvert}{\llVert B\times C\rrVert\llVert C\times A\rrVert},
 \end{eqnarray*}
 and hence
 \[\frac{\sinh(\gamma)}{\sinh(c)}=-\frac{\lvert\det(A,B,C)\rvert}{\llVert A\times B\rrVert\llVert B\times C\rrVert\llVert C\times A\rrVert}.\]
 The same expression arises for $\frac{\sinh(\beta)}{\sinh(b)}$.
 For the third angle, however, we get
 \begin{eqnarray*}
  \sinh(\alpha)&=&\sinh(\olda')\\
  &=&\llVert\oldB'\times\oldC'\rrVert\\
  &=&\llVert B'\times C'\rrVert,
 \end{eqnarray*}
 which results in the minus sign appearing in the equation to prove.}

\subsection*{Sums of Angles and Lengths}

\rm Reasonable statements can only be obtained for the sum of angles in (antipodal) hyperbolic triangles
and for the sum of lengths in spatiolateral triangles.
All the other quantities that could be considered may be infinitely large or infinitely small.
At the end of this section, we are able to finally prove
the invalidity of the triangle equation for contractible spatiolateral triangles.

\Satz{\label{Winkelsumme_hyp}The sum of angles in an (antipodal) hyperbolic triangle is less than $\pi$.}

\rm The proof can be found in Iversen \cite{Ive}.

\Beweis{WLOG assume that $\alpha\geq\beta$ holds.
 The (antipodal) hyperbolic LCA yields
 \begin{eqnarray*}
  \cos(\alpha)&=&\cosh(\alpha)\sin(\beta)\sin(\gamma)-\cos(\beta)\cos(\gamma)\\
  &>&-\cos(\beta+\gamma)\\
  &=&\cos(\lvert\pi-(\beta+\gamma)\rvert).
 \end{eqnarray*}
 Since $0<\beta+\gamma<2\pi$ holds, we have $\lvert\pi-(\beta+\gamma)\rvert\in[0,\pi)$,
 where the cosine function is strictly monotonically decreasing.
 Thus we have
 \[\alpha<\lvert\pi-(\beta+\gamma)\rvert.\]
 The assumption $\alpha\geq\beta$ assures that $\pi-(\beta+\gamma)$ is positive,
 whence we get the desired result.}

\Satz{\label{Satz:gross_klein} The lengths of the sides in a non-contractible spatiolateral triangle sum up
 to a number greater than $2\pi$.
 For contractible spatiolateral triangles, the sum of the lengths of its sides is less than $2\pi$.}

\Beweis{Let firstly the spatiolateral triangle $\{A,B,C\}$ be non-contractible.
 Then the polar triangle $\{A',B',C'\}$ is hyperbolic or antipodal hyperbolic.
 The relation between the angles of the polar triangle and the sides of the original triangle
 which we got by lemma \ref{Satz:H2_S11} tells us
 \[\alpha'=\pi-a,\quad\beta'=\pi-b,\quad\gamma'=\pi-c,\]
 which together with the previous theorem proves our assertion.

 Regarding a contractible spatiolateral triangle, we construct the triangle $\{\oldA',\oldB',\oldC'\}$
 according to theorem \ref{klein_Winkelcos}.
 This theorem tells us that
 $$\cos(a)=-\cos(\alpha'),$$
 $$\cos(b)=\cos(\beta'), \makebox[0pt][l]{ and}$$
 $$\cos(c)=\cos(\gamma')$$
 hold, where $\alpha',\beta',\gamma'$ are the angles of $\{\oldA',\oldB',\oldC'\}$.
 Since the possible values for the sides and angles are in the interval $[0,2\pi]$, we have
 $$\alpha'=\pi-a,\quad\beta'=b,\mbox{ and } \gamma'=c.$$
 Together with theorem \ref{Winkelsumme_hyp}, we get
 \begin{eqnarray*}
  \pi-a+b+c&<&\pi\\
  b+c&<&a<\pi\\
  a+(b+c)&<&2\pi.
 \end{eqnarray*}

 \vspace{-.3cm}}

\Kor{\label{Folg:gross_klein}For non-degenerate contractible spatiolateral de Sitter triangles,
 the triangle inequality does not hold.}

\Beweis{This can be read in the penultimate line of the previous proof.}

\end{document}